\documentclass[11 pt,english]{article}
\usepackage[activeacute]{babel}
\usepackage{amsfonts} 
\usepackage{mathrsfs} 
\usepackage{amssymb, amsmath, amsfonts}
\usepackage{makeidx}
\usepackage{epsfig}
\usepackage{color}
\usepackage[T1]{fontenc}
\usepackage{graphics,graphicx,graphpap}
\usepackage[ansinew]{inputenc}
\usepackage{amsthm}

\newcommand\blfootnote[1]{%
  \begingroup
  \renewcommand\thefootnote{}\footnote{#1}%
  \addtocounter{footnote}{-1}%
  \endgroup
}

\usepackage{tikz}
\usepackage{subfigure}
\usetikzlibrary{calc, positioning, through, math, arrows, backgrounds}
\tikzstyle{vtx}=[inner sep=1pt,draw, shape=circle, font=\tiny]
\tikzstyle{line}=[inner sep=3pt,draw, shape=rectangle, line width = 3pt]
\tikzset{>=stealth}
\tikzstyle{lbl}=[inner sep = 1 pt, fill = white, midway]

\if@twoside \oddsidemargin 6pt \evensidemargin 6pt \marginparwidth 20pt
 \else \oddsidemargin 6pt \evensidemargin 6pt \marginparwidth 20pt \fi
\newcommand{\ZZ}{\mathbb{Z}}

\newcommand{\Aut}{\mathrm{Aut}}
\newcommand{\G}{\Gamma}

\newcommand{\cG}{\mathcal{G}}
\newcommand{\cR}{\mathcal{R}}
\newcommand{\cC}{\mathcal{C}}

\newcommand{\cI}{\mathcal{I}}

\newcommand{\cXa}{\mathcal{X}_a}
\newcommand{\cXb}{\mathcal{X}_b}

\newcommand{\la}{\langle}
\newcommand{\ra}{\rangle}

\newcommand{\No}[1]{{#1}^*}
\newcommand{\vv}[1]{u_{#1}}

\newtheorem{theorem}{Theorem}[section]
\newtheorem{proposition}[theorem]{Proposition}

\newtheorem{lemma}[theorem]{Lemma}
\newtheorem{problem}[theorem]{Problem}

\theoremstyle{definition}

\newtheorem{construction}[theorem]{Construction}

\begin{document}

\begin{center}
\Large{\textbf{Cubic factor-invariant graphs of bialternating cycle quotient type}} \\ [+4ex]
\end{center}

\begin{center}
Primo\v z \v Sparl
\\

\medskip
{\it {\small
University of Ljubljana, Faculty of Education, Ljubljana, Slovenia\\
University of Primorska, Institute Andrej Maru\v si\v c, Koper, Slovenia\\
Institute of Mathematics, Physics and Mechanics, Ljubljana, Slovenia\\
}}
\end{center}

\blfootnote{
Email address: 
primoz.sparl@pef.uni-lj.si
}


\hrule

\begin{abstract}
In 2019, investigation of the so-called factor-invariant cubic graphs was initiated by Alspach, Khodadadpour and Kreher. For a cubic graph $\G$ and a vertex-transitive subgroup $G$ of $\Aut(\G)$, a $2$-factor $\cC$ of $\G$ is said to be {\em $G$-invariant} if the set $\cC$ is preserved by each element of $G$. Investigations of factor-invariant cubic graphs therefore contribute to the rapidly growing theory on cubic vertex-transitive graphs, providing a better insight into the structure of such graphs.  

Initially, the examples where $\cC$ consists of a single or just two cycles were analyzed. In a recent paper by Brian Alspach and the author of this paper, the investigation of the examples for which the corresponding quotient graph $\G_\cC$ of $\G$ with respect to $\cC$ is a cycle was initiated. Moreover, the graphs of the so-called {\em alternating cycle quotient type} were classified. 

In this paper, the remaining examples, that is the graphs of the {\em bialternating cycle quotient type}, are classified. It is shown that they belong to a previously unknown infinite $5$-parametric family of graphs of girth at most $10$ and that they are Cayley graphs of groups with respect to three involutions. 
\end{abstract}
\hrule

\begin{quotation}
\noindent {\em \small Keywords: cubic; vertex-transitive; factor-invariant}
\end{quotation}

\section{Introduction}
\label{sec:Intro}

This paper is a contribution to the rapidly growing theory on cubic vertex-transitive graphs (see~\cite{AlsSpa24, BarGraSpi25, DobHujImrOrt25, GuoMoh24, KutMarMikSpa24, LiKwoZho23, PotTol23} for some of the most recent results). The obtained results complement those of~\cite{AlsSpa24}, where cubic factor-invariant graphs of alternating cycle quotient type were classified -- see the next paragraph for the terminology. In this way a new step in the efforts concerning~\cite[Problem~1.1]{AlsSpa24} is made. This problem, which originates from a more specific question of Bojan Mohar to Brian Alspach from 2016 (see~\cite{AlsKreKho19}), asks for a structural characterization of cubic vertex-transitive graphs admitting a partition of their edge-sets into a $2$-factor and a $1$-factor which is preserved by some vertex-transitive subgroup of the automorphism group. What is meant here is a description of such graphs in terms of suitable parameter sets, where the parameters determine the number and length of the cycles from the corresponding $2$-factor $\cC$, as well as describe how different cycles from $\cC$ are linked among themselves. In this way, interesting new families of cubic vertex-transitive graphs can be discovered. 

Besides the two most basic possibilities, where the corresponding $2$-factor $\cC$ consists of just one or two cycles (for a complete classification, see~\cite{AlsKreKho19} and~\cite{AlsDobKhoSpa22}, respectively), one of the most natural possibilities for the above mentioned structure of such a cubic graph is that each member of $\cC$ is adjacent to just two other members of $\cC$. In such a case the graph is said to be of {\em cycle quotient type} with respect to $\cC$ (see Section~\ref{sec:Prelim} for a precise definition). In~\cite{AlsSpa24}, computational data based on~\cite{PotSpiVer13} was presented, which suggests that investigating such graphs makes sense. It was proved in~\cite{AlsSpa24} that the graphs of cycle quotient type come in two essentially different ``subtypes'', namely they are either of the so-called {\em alternating} or of {\em bialternating} cycle quotient type. Those of alternating cycle quotient type were thoroughly investigated in~\cite{AlsSpa24}, where it was proved that all such graphs belong to a certain $4$-parametric family of graphs $\cXa(m,n,k,\ell)$. It was also pointed out that this new family is a very natural generalization of the widely studied family of generalized Petersen graphs (as well as of the family of Honeycomb toroidal graphs~\cite{Als21, Spa22}). In addition, the automorphism group of each member $\cXa(m,n,k,\ell)$ of this family was determined.

The aim of this paper is a thorough investigation of the examples of the bialternating cycle quotient type. Our main result, Theorem~\ref{the:main}, gives a complete characterization of such graphs. Besides the well-known prism graphs and the M\"obius ladders, all such graphs either belong to two very specific families of examples (denoted $\cXb^1(m)$ and $\cXb^2(m)$) having cycles from $\cC$ of length $8$, or belong to a new $5$-parametric family of graphs $\cXb(m,n,a,b,\ell)$, which to the best of our knowledge was not known and/or studied thus far. We show that all such graphs are of girth at most $10$ and are Cayley graphs with respect to a connection set consisting of three involutions. For some we also determine their full automorphism groups.

\section{Preliminaries}
\label{sec:Prelim}

For a comprehensive introduction of the relevant theory and notation we refer the reader to~\cite{AlsSpa24}, but for self-completeness we repeat the definitions and notational conventions from~\cite{AlsSpa24} that we will be using throughout this paper.

For a graph $\G$ we denote its vertex set and edge set by $V(\G)$ and $E(\G)$, respectively, where we sometimes simply write $V$ and $E$ if the graph $\G$ is clear from the context. For a subgroup $G$ of the automorphism group $\Aut(\G)$ of a graph $\G$, the graph $\G$ is said to be {\em $G$-vertex-transitive} or {\em $G$-edge-transitive}, if the natural action of $G$ on the set $V(\G)$ or $E(\G)$, respectively, is transitive. Similarly, for an integer $s \geq 1$, the graph $\G$ is said to be {\em $G$-$s$-arc-transitive}, if the natural action of $G$ on the set of all sequences of vertices $(v_0,v_1,\ldots , v_s)$ of $\G$, for which any two consecutive vertices are adjacent and any three consecutive vertices are pairwise distinct (such sequences are called {\em $s$-arcs}), is transitive. Moreover, if the action of $G$ on the set of all $s$-arcs is regular, we say that $\G$ is {\em $G$-$s$-arc-regular}. In the case that $G = \Aut(\G)$, we omit the prefix $\Aut(\G)$ and simply say that the graph $\G$ is {\em vertex-transitive}, {\em edge-transitive}, {\em $s$-arc-transitive} or {\em $s$-arc-regular}. We also abbreviate the terms $1$-arc-transitive and $1$-arc-regular to arc-transitive and arc-regular, respectively.

Throughout the paper the Greek letters such as $\alpha, \beta, \gamma, \varphi$ and $\rho$ will denote automorphisms of graphs while $\delta$ and $\varepsilon$ will usually be integers from $\{0,1\}$ or $\{-1,1\}$. This should cause no confusion.

We will constantly be working with elements of the residue class ring $\ZZ_n$ for some positive integer $n$, as well as with integers. Often we will write things such as $t\ell + 2 = 0$, where $\ell$ will be an element of $\ZZ_n$ and $t$ will be an integer. By this we mean that the corresponding equality holds in $\ZZ_n$, or in other words that viewing $\ell$ as an integer (any one in the corresponding residue class), $t\ell + 2$ is divisible by $n$. We also make the convention that for an integer (or an element of $\ZZ_n$), say $k$, we abbreviate things like $k \in \{-2,2\}$ to $k = \pm 2$ and $k \notin \{-2,2\}$ to $k \neq \pm 2$.
\medskip

In all subsequent sections $\G$ will usually denote a cubic vertex-transitive graph admitting a vertex-transitive group $G \leq \Aut(\G)$ and a partition of the edge set of $\G$ into a $2$-factor $\cC$ and a $1$-factor $\cI$ that is preserved by $G$. In this case we say that this partition (as well as $\cC$) is {\em $G$-invariant}. The corresponding {\em quotient graph} $\G_\cC$ of $\G$ {\em with respect to} $\cC$ is then the graph with vertex set $\cC$ and in which two distinct $C, C' \in \cC$ are adjacent whenever there are vertices $u$ of $C$ and $v$ of $C'$ which are adjacent in $\G$. As explained in the Introduction, we will be working with examples of {\em cycle quotient type} (with respect to $\cC$) which are the ones for which $|\cC| \geq 3$ and $\G_\cC$ is a cycle. 

Throughout the paper we let $n$ be the length of the cycles from $\cC$ and we let $m = |\cC|$. This of course means that $\G$ is of order $mn$. We denote the cycles from $\cC$ by $C_i$, $i \in \ZZ_m$, in such a way that $C_i \sim C_{i-1}, C_{i+1}$ in $\G_\cC$ (indices computed modulo $m$). Moreover, for each $i \in \ZZ_m$ we let $V_i = \{\vv{i,j} \colon j \in \ZZ_n\}$ be the vertex set of $C_i$ (and we therefore always compute the indices of $\vv{i,j}$ modulo $m$ in the first and modulo $n$ in the second component). Since $\G$ is cubic, each vertex $v$ of $\G$ has precisely one neighbor which is in a different member of $\cC$ than $v$. We denote this vertex by $\No{v}$ and call it the {\em outside neighbor} of $v$. Note that $\No{(\No{v})} = v$. The following straightforward observation will be useful.

\begin{lemma}
\label{le:aut_star}
Let $\G$ be a cubic vertex-transitive graph admitting a vertex-transitive group $G \leq \Aut(\G)$ and a $G$-invariant partition of the edge set of $\G$ into a $2$-factor with at least two components and a $1$-factor. Then $\varphi(\No{v}) = \No{(\varphi(v))}$ holds for each $v \in V(\G)$ and each $\varphi \in G$.
\end{lemma}

\section{The graphs $\cXb(m,n,a,b,\ell)$}
\label{sec:graphs}

Throughout this section let $\G$ be a connected cubic vertex-transitive graph admitting a vertex-transitive group $G \leq \Aut(\G)$ and a $G$-invariant partition of its edge set into a $2$-factor $\cC$ and a $1$-factor $\cI$ such that $|\cC| \geq 3$ and $\G$ is of cycle quotient type with respect to $\cC$. 

Suppose the vertices of $\G$ have been labeled as described in the previous section. With no loss of generality we can assume that $\vv{0,j} \sim \vv{0,j+1}$ for all $j \in \ZZ_n$ and that $\No{\vv{0,0}} = \vv{1,0}$. As explained in the Introduction, the examples for which $\No{\vv{0,1}}, \No{\vv{0,-1}} \in V_{m-1}$ have been classified in~\cite{AlsSpa24}, and so we assume that at least one of $\No{\vv{0,1}}, \No{\vv{0,-1}}$ is in $V_1$. Since $\G$ is $G$-vertex-transitive and $\cC$ is $G$-invariant, the restriction of the setwise stabilizer of $V_0$ in $G$ to $V_0$ is transitive, and so is isomorphic to a transitive subgroup of the dihedral group $D_n$ of order $2n$. There thus exists some $\rho \in G$ preserving $C_0$, whose action on $V_0$ is the $2$-step rotation given by 
\begin{equation}
\label{eq:rho_V_0}
	\rho(\vv{0,j}) = \vv{0,j+2}\ \text{for\ all}\ j \in \ZZ_n.
\end{equation}
Since $|\cC| \geq 3$, this clearly implies that precisely one of $\vv{0,-1}, \vv{0,1}$ has its outside neighbour in $V_{1}$. With no loss of generality we assume that $\No{\vv{0,-1}} \in V_{m-1}$ and $\No{\vv{0,1}} \in V_1$. As mentioned in the Introduction, we say that $\G$ is of {\em bialternating cycle quotient type} with respect to $\cC$ in this case. Unlike in the case of the graphs $\cXa(m,n,k,\ell)$ of alternating cycle quotient type from~\cite{AlsSpa24}, we this time assume that the vertices of $\G$ have been labeled in such a way that 
\begin{equation}
\label{eq:non-links}
	\vv{i,j} \sim \vv{i,j-1},\vv{i,j+1}\ \text{for\ all}\ i \in \ZZ_m,\ j \in \ZZ_n, 
\end{equation}
which we can do with no loss of generality since each $C_i$ is a cycle of length $n$.
The following observation will be very useful.

\begin{lemma}
\label{le:G_regular}
With the assumptions and notation of this section, the length $n$ of the cycles from $\cC$ is divisible by $4$ and the group $G$ acts regularly on the vertex set of $\G$. Moreover, for each $i \in \ZZ_m$, the setwise stabilizer of $V_i$ in $G$ is isomorphic to the dihedral group $D_{n/2}$ of order $n$.
\end{lemma}

\begin{proof}
Suppose $\varphi \in G$ fixes the vertex $\vv{0,0}$. By Lemma~\ref{le:aut_star}, it then also fixes its outside neighbour $\vv{1,0}$, and thus preserves each of $V_0$ and $V_1$. Consequently, it also preserves each $V_i$, $i \in \ZZ_m$. Since $\No{\vv{0,1}} \in V_1$, while $\No{\vv{0,-1}} \notin V_1$, the vertices $\vv{0,1}$ and $\vv{0,-1}$ are also fixed by $\varphi$. This shows that $\varphi$ fixes each neighbour of $\vv{0,0}$. Since $G$ is vertex-transitive and connected, an inductive argument shows that $\varphi$ is the identity, thus proving that $G$ indeed acts regularly on the vertex set of $\G$.

To complete the proof let $\rho$ be the (unique) element of $G$ satisfying~\eqref{eq:rho_V_0}. Since $\No{\vv{0,-1}} \in V_{m-1}$ and $\No{\vv{0,1}} \in V_1$, the induced action of $\rho$ on $\G_\cC$ fixes the vertex $C_0$ and interchanges $C_1$ with $C_{m-1}$. It is thus clear that $\rho$ is of even order, implying that $4 \mid n$, as claimed. Finally, since the setwise stabilizer of $V_0$ in $G$ is transitive on the cycle $C_0$ but $G$ does not contain a $1$-step rotation of $C_0$ (since $\No{\vv{0,0}}, \No{\vv{0,1}} \in V_1$ but $\No{\vv{0,-1}} \notin V_1$), the restriction of this stabilizer to $V_0$ must be the dihedral group $D_{n/2}$ (acting regularly on $V_0$). The same then clearly holds for the setwise stabilizer of each $V_i$ in $G$. Since $G$ is regular, this completes the proof.
\end{proof} 

Applying $\rho^2$ from the above proof (note that $\rho^2 = 1$ in the case of $n = 4$) also yields the following fact, which we record for future reference: 
\begin{align}
	\No{\vv{0,j}} \in V_1 & \iff j \equiv 0 \pmod{4}\ \text{or}\ j \equiv 1 \pmod{4},\label{eq:C0n1}\\
	\No{\vv{0,j}} \in V_{m-1} & \iff j \equiv 2 \pmod{4}\ \text{or}\ j \equiv 3 \pmod{4}. \label{eq:C0n2}
\end{align}

\subsection{The automorphisms $\rho$, $\alpha$ and $\gamma$}
\label{subsec:agr}

We now consider three specific elements of the regular group $G$, which will enable us to completely determine the structure of the graph $\G$. By Lemma~\ref{le:G_regular}, there is a unique $\rho \in G$ satisfying~\eqref{eq:rho_V_0}. We thus make an agreement that whenever we speak of `the $\rho$' from $G$ we are referring to this unique element. Similarly, there exists a unique $\alpha \in G$ mapping $\vv{0,0}$ to its outside neighbour $\vv{1,0}$. By Lemma~\ref{le:aut_star} it then maps $\vv{1,0}$ back to $\vv{0,0}$, and is thus an involution by Lemma~\ref{le:G_regular}. Recall our assumption that~\eqref{eq:non-links} holds. Since precisely one of $\No{\vv{1,1}}$ and $\No{\vv{1,-1}}$ is in $V_0$, we can assume that the vertices of $C_1$ have been labeled in such a way that $\No{\vv{1,1}} \in V_0$. In other words, we can assume that $\alpha$ interchanges $\vv{0,1}$ with $\vv{1,1}$, which in turn implies that  
\begin{equation}
\label{eq:alpha}
	\alpha(\vv{0,j}) = \vv{1,j}\ \text{and}\ \alpha(\vv{1,j}) = \vv{0,j}\ \text{for\ all}\ j \in \ZZ_n.
\end{equation} 
We again speak of `the $\alpha$' from $G$ when we are referring to the unique $\alpha \in G$ satisfying~\eqref{eq:alpha}.

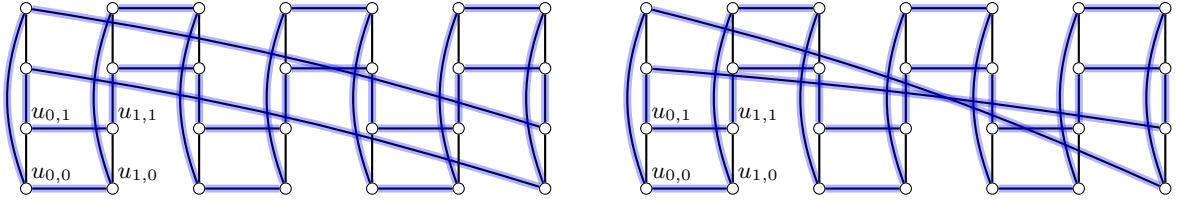
\begin{figure}[h]
\begin{center}
\subfigure
{
\begin{tikzpicture}[scale = .5]
\foreach \i in {0,1,2,3,4,5,6}{
\foreach \j in {0,1,2,3}{
\node[vtx, fill = white, inner sep = 1.5pt,] (A\i\j) at (2.3*\i, 1.6*\j) {};
}}
\node[above right = -5pt of A00] {\footnotesize $\vv{0,0}$};
\node[above right = -5pt of A10] {\footnotesize $\vv{1,0}$};
\node[above right = -5pt of A01] {\footnotesize $\vv{0,1}$};
\node[above right = -5pt of A11] {\footnotesize $\vv{1,1}$};
\begin{scope}[on background layer]
\foreach \j in {0,1,2}{
\foreach \i in {0,1,2,3,4,5,6}{
\draw[thick] let \n1 = {int(mod(\j+1, 4))} in (A\i\j) -- (A\i\n1);
}}
\foreach \i in {0,1,2,3,4,5,6}{
\draw[thick] (A\i 0) to[bend left = 20] (A\i 3);
}
\foreach \i in {0,2,4}{
\foreach \j in {0,1}{
\draw[thick] let \n1 = {int(mod(\i+1, 7))} in (A\i\j) -- (A\n1\j);
\draw[thick] let \n1 = {int(mod(\i+1, 7))}, \n2 = {int(mod(\i+2, 7))}, \n3 = {int(mod(\j + 2,4))} in (A\n1\n3) -- (A\n2\n3);
}}
\draw[thick] (A60) to[bend right = 4] (A02);
\draw[thick] (A61) to[bend right = 4] (A03);
\foreach \i in {0,2,4}{
\draw[opacity = .3, blue, line width = 3 pt] let \n1 = {int(mod(\i+1,7))} in (A\i0) -- (A\n10) to[bend left = 20] (A\n13);
\draw[opacity = .3, blue, line width = 3 pt] let \n1 = {int(mod(\i+1,7))}, \n2 = {int(mod(\i+2,7))} in (A\n13) -- (A\n23) to[bend right = 20] (A\n20);
\draw[opacity = .3, blue, line width = 3 pt] let \n1 = {int(mod(\i+1,7))}, \n2 = {int(mod(\i+2,7))} in (A\i1) -- (A\n11) -- (A\n12) -- (A\n22) -- (A\n21);
}
\draw[opacity = .3, blue, line width = 3 pt] (A60) to[bend right = 4] (A02);
\draw[opacity = .3, blue, line width = 3 pt] (A61) to[bend right = 4] (A03);
\draw[opacity = .3, blue, line width = 3 pt] (A00) to[bend left = 20] (A03);
\draw[opacity = .3, blue, line width = 3 pt] (A01) -- (A02);
\end{scope}
\end{tikzpicture}
}
\hskip 5mm
\subfigure
{
\begin{tikzpicture}[scale = .5]
\foreach \i in {0,1,2,3,4,5,6}{
\foreach \j in {0,1,2,3}{
\node[vtx, fill = white, inner sep = 1.5pt,] (A\i\j) at (2.3*\i, 1.6*\j) {};
}}
\node[above right = -5pt of A00] {\footnotesize $\vv{0,0}$};
\node[above right = -5pt of A10] {\footnotesize $\vv{1,0}$};
\node[above right = -5pt of A01] {\footnotesize $\vv{0,1}$};
\node[above right = -5pt of A11] {\footnotesize $\vv{1,1}$};
\begin{scope}[on background layer]
\foreach \j in {0,1,2}{
\foreach \i in {0,1,2,3,4,5,6}{
\draw[thick] let \n1 = {int(mod(\j+1, 4))} in (A\i\j) -- (A\i\n1);
}}
\foreach \i in {0,1,2,3,4,5,6}{
\draw[thick] (A\i 0) to[bend left = 20] (A\i 3);
}
\foreach \i in {0,2,4}{
\foreach \j in {0,1}{
\draw[thick] let \n1 = {int(mod(\i+1, 7))} in (A\i\j) -- (A\n1\j);
\draw[thick] let \n1 = {int(mod(\i+1, 7))}, \n2 = {int(mod(\i+2, 7))}, \n3 = {int(mod(\j + 2,4))} in (A\n1\n3) -- (A\n2\n3);
}}
\draw[thick] (A60) to[bend right = 5] (A03);
\draw[thick] (A61) to[bend right = 2] (A02);
\foreach \i in {0,2,4}{
\draw[opacity = .3, blue, line width = 3 pt] let \n1 = {int(mod(\i+1,7))} in (A\i0) -- (A\n10) to[bend left = 20] (A\n13);
\draw[opacity = .3, blue, line width = 3 pt] let \n1 = {int(mod(\i+1,7))}, \n2 = {int(mod(\i+2,7))} in (A\n13) -- (A\n23) to[bend right = 20] (A\n20);
\draw[opacity = .3, blue, line width = 3 pt] let \n1 = {int(mod(\i+1,7))}, \n2 = {int(mod(\i+2,7))} in (A\i1) -- (A\n11) -- (A\n12) -- (A\n22) -- (A\n21);
}
\draw[opacity = .3, blue, line width = 3 pt] (A60) to[bend right = 5] (A03);
\draw[opacity = .3, blue, line width = 3 pt] (A61) to[bend right = 2] (A02);
\draw[opacity = .3, blue, line width = 3 pt] (A00) to[bend left = 20] (A03);
\draw[opacity = .3, blue, line width = 3 pt] (A01) -- (A02);
\end{scope}
\end{tikzpicture}
}
\caption{The M\"obius ladder and the prism resulting from the case $m = 7$ and $n = 4$.}
\label{fig:n=4}
\end{center}
\end{figure}

To describe the edges from the $1$-factor $\cI$, which we refer to as {\em links} throughout the paper, let $a \in \ZZ_n$ be the unique element such that $\No{\vv{0,1}} = \vv{1,a}$. Lemma~\ref{le:aut_star} and \eqref{eq:alpha} then imply that $\No{\vv{0,a}} = \vv{1,1}$, and so \eqref{eq:C0n1} implies that one of $a$ and $a-1$ is divisible by $4$. If $n = 4$, this forces $a = 1$. Using the fact that $G$ is vertex-transitive and preserves $\cC$, it is easy to see that in this case $\G$ is either the M\"obius ladder or the prism of order $4m$ (see Figure~\ref{fig:n=4} for the two examples with $m = 7$, where the edges of the graph that do not belong to two $4$-cycles of the graph have been highlighted). Note that in both of these cases the full automorphism group $\Aut(\G)$ has vertex stabilizers of order $2$ and does not preserve $\cC$. However, it can be verified that $\Aut(\G)$ indeed contains a regular subgroup $G$ preserving $\cC$.
%

From now on we will therefore assume that $n \geq 8$. By~\eqref{eq:C0n1}, there exists an element $b \in \ZZ_n \setminus \{a\}$ such that $\No{\vv{0,4}} = \vv{1,b}$. Lemma~\ref{le:aut_star} and~\eqref{eq:alpha} then imply that $\No{\vv{0,b}} = \vv{1,4}$, and so one of $b$ and $b-1$ is divisible by $4$. Note also that~\eqref{eq:C0n1} and~\eqref{eq:alpha} force $\No{\vv{1,2}}, \No{\vv{1,3}} \in V_2$. 

Let $\gamma$ be the unique element of $G$ mapping $\vv{0,0}$ to $\vv{1,2}$. Since $\No{\vv{0,0}} = \vv{1,0} \in V_1$ and $\No{\vv{1,2}} \in V_2$, it follows that $\gamma$ maps each $C_i$ to $C_{i+1}$ (where $i+1$ is computed modulo $m$). Since $\No{\vv{0,1}} \in V_1$ and $\No{\vv{1,3}} \in V_2$ while $\No{\vv{1,1}} \in V_0$, we deduce that $\gamma(\vv{0,1}) = \vv{1,3}$, and consequently $\gamma(\vv{0,j}) = \vv{1,j+2}$ for all $j \in \ZZ_n$. We can assume that the vertices of $C_2$ have been labeled in such a way that $\No{\vv{1,2}} = \vv{2,2}$ and $\No{\vv{2,3}} \in V_1$. Then Lemma~\ref{le:aut_star} implies that $\gamma(\vv{1,0}) = \vv{2,2}$, and since $\No{\vv{1,1}} \in V_0$, we also find that $\gamma(\vv{1,1}) = \vv{2,3}$. Thus $\gamma(\vv{1,j}) = \vv{2,j+2}$ for all $j \in \ZZ_n$. Consequently, Lemma~\ref{le:aut_star} implies that $\No{\vv{1,3}} = \gamma(\vv{1,a}) = \vv{2,2+a}$ and $\No{\vv{1,6}} = \gamma(\vv{1,b}) = \vv{2,2+b}$. Continuing in this way we can thus assume that the vertices of $\G$ have been labeled in such a way that 
\begin{equation}
\label{eq:gamma}
\gamma(\vv{i,j}) = \vv{i+1, j+2}\ \text{for\ all}\ i \in \{0,1,\ldots , m-2\}\ \text{and}\ j \in \ZZ_n,
\end{equation}
and that 
\begin{equation}
\label{eq:links_init}
	\No{\vv{i,2i}} = \vv{i+1,2i},\ \No{\vv{i,2i+1}} = \vv{i+1,2i+a}\ \text{and}\ \No{\vv{i,2i+4}} = \vv{i+1,2i+b}
\end{equation}
holds for each $i \in \{0,1,\ldots , m-2\}$. Recall that $\gamma$ maps $V_{m-1}$ to $V_0$ and note that~\eqref{eq:gamma} and~\eqref{eq:links_init} imply that $\No{\vv{m-1,2(m-1)}}, \No{\vv{m-1,2(m-1)+1}} \in V_0$. Let $\ell \in \ZZ_n$ be such that $\No{\vv{m-1,2(m-1)}} = \vv{0,\ell}$ and note that by \eqref{eq:C0n2}, one of $\ell-2$ and $\ell-3$ is divisible by $4$. Since $\No{\vv{m-2,2(m-2)}} = \vv{m-1,2(m-2)}$, it thus follows that $\gamma(\vv{m-1,2(m-2)}) = \vv{0,\ell}$. Therefore, if $\ell$ is even (in which case $\No{\vv{0,\ell+1}} \in V_{m-1})$, then
\begin{equation}
\label{eq:gamma_m1}
	\gamma(\vv{m-1,j}) = \vv{0,\ell + j - 2(m-2)}\ \text{for\ all}\ j \in \ZZ_n,
\end{equation}
while if $\ell$ is odd (in which case $\No{\vv{0,\ell-1}} \in V_{m-1})$, then
\begin{equation}
\label{eq:gamma_m2}
	\gamma(\vv{m-1,j}) = \vv{0,\ell - j + 2(m-2)}\ \text{for\ all}\ j \in \ZZ_n.
\end{equation}

\subsection{The four types for the $n = 8$ case}
\label{subsec:n=8}

It turns out that the situation can be somewhat specific when $n = 8$, which is why this possibility needs to be analyzed with particular care. We defer the complete analysis of the $n = 8$ case to Section~\ref{sec:n=8}, but since two of the four subcases can be consider along the general $n > 8$ case, we make some remarks, which will help us to do so. 

\begin{figure}[h]
\begin{center}
\subfigure
{
\begin{tikzpicture}[scale = .6]
\foreach \i in {0,1,2}{
\foreach \j in {0,1,2,3,4,5,6,7}{
\node[vtx, fill = white, inner sep = 1.5pt,] (A\i\j) at (2*\i, 1*\j) {};
}}
\node[below = 3pt of A10] {\footnotesize $a = 1$, $b = 4$};
\begin{scope}[on background layer]
\foreach \j in {0,1,2,3,4,5,6}{
\foreach \i in {0,1,2}{
\draw[thick] let \n1 = {int(mod(\j+1, 8))} in (A\i\j) -- (A\i\n1);
}}
\foreach \i in {0,1,2}{
\draw[thick] (A\i 0) to[bend left = 15] (A\i 7);
}
\foreach \j in {0,1}{
\draw[thick] (A0\j) -- (A1\j);
\draw[thick] let \n1 = {int(mod(\j+4, 8))} in (A0\n1) -- (A1\n1);
\draw[thick] let \n1 = {int(mod(\j+2, 8))} in (A1\n1) -- (A2\n1);
\draw[thick] let \n1 = {int(mod(\j+6, 8))} in (A1\n1) -- (A2\n1);
}
\end{scope}
\end{tikzpicture}
}
\hskip 5mm
\subfigure
{
\begin{tikzpicture}[scale = .6]
\foreach \i in {0,1,2}{
\foreach \j in {0,1,2,3,4,5,6,7}{
\node[vtx, fill = white, inner sep = 1.5pt,] (A\i\j) at (2*\i, 1*\j) {};
}}
\node[below = 3pt of A10] {\footnotesize $a = 5$, $b = 4$};
\begin{scope}[on background layer]
\foreach \j in {0,1,2,3,4,5,6}{
\foreach \i in {0,1,2}{
\draw[thick] let \n1 = {int(mod(\j+1, 8))} in (A\i\j) -- (A\i\n1);
}}
\foreach \i in {0,1,2}{
\draw[thick] (A\i 0) to[bend left = 15] (A\i 7);
}
\draw[thick] (A00) -- (A10);
\draw[thick] (A12) -- (A22);
\draw[thick] (A01) -- (A15);
\draw[thick] (A04) -- (A14);
\draw[thick] (A05) -- (A11);
\draw[thick] (A16) -- (A26);
\draw[thick] (A13) -- (A27);
\draw[thick] (A17) -- (A23);
\end{scope}
\end{tikzpicture}
}
\hskip 5mm
\subfigure
{
\begin{tikzpicture}[scale = .6]
\foreach \i in {0,1,2}{
\foreach \j in {0,1,2,3,4,5,6,7}{
\node[vtx, fill = white, inner sep = 1.5pt,] (A\i\j) at (2*\i, 1*\j) {};
}}
\node[below = 3pt of A10] {\footnotesize $a = 4$, $b = 1$};
\begin{scope}[on background layer]
\foreach \j in {0,1,2,3,4,5,6}{
\foreach \i in {0,1,2}{
\draw[thick] let \n1 = {int(mod(\j+1, 8))} in (A\i\j) -- (A\i\n1);
}}
\foreach \i in {0,1,2}{
\draw[thick] (A\i 0) to[bend left = 15] (A\i 7);
}
\draw[thick] (A00) -- (A10);
\draw[thick] (A12) -- (A22);
\draw[thick] (A01) -- (A14);
\draw[thick] (A04) -- (A11);
\draw[thick] (A05) -- (A15);
\draw[thick] (A16) -- (A23);
\draw[thick] (A13) -- (A26);
\draw[thick] (A17) -- (A27);
\end{scope}
\end{tikzpicture}
}
\hskip 5mm
\subfigure
{
\begin{tikzpicture}[scale = .6]
\foreach \i in {0,1,2}{
\foreach \j in {0,1,2,3,4,5,6,7}{
\node[vtx, fill = white, inner sep = 1.5pt,] (A\i\j) at (2*\i, 1*\j) {};
}}
\begin{scope}[on background layer]
\foreach \j in {0,1,2,3,4,5,6}{
\foreach \i in {0,1,2}{
\draw[thick] let \n1 = {int(mod(\j+1, 8))} in (A\i\j) -- (A\i\n1);
}}
\node[below = 3pt of A10] {\footnotesize $a = 1$, $b = 5$};
\foreach \i in {0,1,2}{
\draw[thick] (A\i 0) to[bend left = 15] (A\i 7);
}
\draw[thick] (A00) -- (A10);
\draw[thick] (A12) -- (A22);
\draw[thick] (A01) -- (A11);
\draw[thick] (A04) -- (A15);
\draw[thick] (A05) -- (A14);
\draw[thick] (A16) -- (A27);
\draw[thick] (A13) -- (A23);
\draw[thick] (A17) -- (A26);
\end{scope}
\end{tikzpicture}
}
\caption{The four possibilities for the pair $(a,b)$ in the case of $n = 8$.}
\label{fig:n=8nalt}
\end{center}
\end{figure}
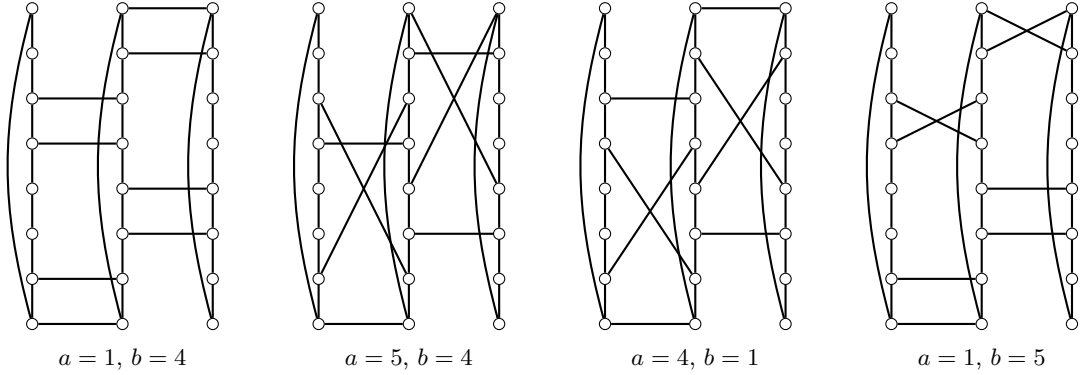

Recall that one of $a$ and $a-1$ and one of $b$ and $b-1$ is divisible by $4$. Because of~\eqref{eq:alpha}, there are thus four essentially different possibilities for the pair $(a,b)$ when $n = 8$, namely $(a,b) \in \{(1,4), (5,4), (4,1), (1,5)\}$. To illustrate these possibilities and our assumption on the labeling of the vertices of $\G$ in general we present the subgraph induced on $V_0 \cup V_1 \cup V_2$ (without indicating the edges between $V_0$ and $V_2$ should $m = 3$ hold) for each of the four possibilities for the pair $(a,b)$ in Figure~\ref{fig:n=8nalt}. Observe that the subgraph of $\G$ induced on $V_0 \cup V_1$ is bipartite if and only if $b = 4$.

Suppose first that $\ell$ is odd. Then \eqref{eq:gamma} and~\eqref{eq:gamma_m2} imply that $\gamma$ maps the pair of adjacent vertices $\vv{m-1,2(m-1)}$ and $\vv{0,\ell}$ to the pair $\vv{0,\ell-2}$ and $\vv{1,\ell+2}$. Since $\ell \in \{3,7\}$, we see that $a = 5$ and $b = 4$ must hold in this case. 

Suppose now that $\ell$ is even. Then $\gamma$ maps the pair of adjacent vertices $\vv{m-1,2(m-1)}$ and $\vv{0,\ell}$ to the pair $\vv{0,\ell+2}$ and $\vv{1,\ell+2}$. If $\ell = 2$, we thus obtain $b = 4$, while if $\ell = 6$, we get no extra condition. We have thus established the following result.

\begin{lemma}
\label{le:n=8}
With the notation from this section, suppose $n = 8$. Then one of the following holds:
\begin{itemize}
\itemsep = 0pt
	\item $\ell$ is odd and $(a,b) = (5,4)$.
	\item $\ell = 2$ and $(a,b) \in \{(1,4), (5,4)\}$.
	\item $\ell = 6$ and $(a,b) \in \{(1,4), (5,4), (4,1), (1,5)\}$.
\end{itemize} 
\end{lemma}

\subsection{The links}

As mentioned in the previous section, the $n = 8$ case needs to be considered with particular care, so assume for the time being that $n > 8$. Recall that $n$ is divisible by $4$, so that we can write $n = 4n_0$ for some $n_0 \geq 3$. By Lemma~\ref{le:G_regular}, $\rho^2$ is of order $n_0$, and is thus not an involution. Since $\No{\vv{0,0}} = \vv{1,0}$ and $\No{\vv{0,4}} = \vv{1,b}$, the restriction of $\rho^2$ to $C_1$ is a rotation of order $n_0$, mapping each $\vv{1,j}$ to $\vv{1,j+b}$. Consequently, $\gcd(n,b) = 4$, and so 
\begin{equation}
\label{eq:b}
	b = 4b_0\ \text{for\ some}\ b_0\ \text{with}\ 1 \leq b_0 < n_0\ \text{and}\ \gcd(n_0, b_0) = 1.
\end{equation}
By Lemma~\ref{le:aut_star}, this also shows that 
\begin{equation}
\label{eq:adjnalt}
	\No{\vv{0,4j_0}} = \vv{1,bj_0} \  \text{and} \  \No{\vv{0,4j_0+1}} = \vv{1,bj_0+a}\ \text{for\ all}\ j_0\ \text{with}\ 0 \leq j_0 < n_0.
\end{equation}
Furthermore,~\eqref{eq:alpha} implies that $\No{\vv{0,b}} = \vv{1,4}$ and $\No{\vv{0,a}} = \vv{1,1}$, and so \eqref{eq:b} and \eqref{eq:adjnalt} yield 
\begin{equation}
\label{eq:b_a_cond}
	bb_0 \equiv 4 \pmod{n},\ a = 4a_0 + 1\ \text{for\ some}\ a_0\ \text{with}\ 0 \leq a_0 < n_0,\ \text{and}\ ba_0 + a \equiv 1 \pmod{n}.
\end{equation} 
That \eqref{eq:b},~\eqref{eq:adjnalt} and~\eqref{eq:b_a_cond} all hold was established based on the assumption that $n > 8$ holds. Note however that even if $n = 8$ and $b = 4$ (in which case $a \in \{1,5\}$), each of \eqref{eq:b},~\eqref{eq:adjnalt} and~\eqref{eq:b_a_cond} still holds. Until Section~\ref{sec:n=8} we will therefore assume that either $n > 8$ or that $n = 8$ with $b = 4$. Observe that in all of these cases the subgraph of $\G$ induced on $V_0 \cup V_1$ is bipartite.

With this assumption in mind, we can now complete the description of the links of $\G$ quite easily. Consider again the automorphism $\gamma \in G$ for which~\eqref{eq:gamma} and one of~\eqref{eq:gamma_m1} and~\eqref{eq:gamma_m2} holds (depending on whether $\ell$ is even or odd, respectively). Using~\eqref{eq:gamma} and~\eqref{eq:adjnalt} we first find that the following holds for each $i$ with $0 \leq i < m-1$, each $j_0$ with $0 \leq j_0 < n_0$, and each $\delta \in \{0,1\}$:
\begin{equation}
\label{eq:links}
	\vv{i,2i+4j_0+\delta} \sim \vv{i+1,2i+bj_0+\delta a}.
\end{equation}
The links between $V_{m-1}$ and $V_0$ depend on the parity of $\ell$. If $\ell$ is even, then~\eqref{eq:gamma_m1} and~\eqref{eq:links} imply that 
\begin{equation}
\label{eq:jumps1}
	\vv{m-1,2(m-1) + 4j_0+\delta} \sim \vv{0,\ell + bj_0+\delta a}
\end{equation}
holds for each $j_0$ with $0 \leq j_0 < n_0$ and $\delta \in \{0,1\}$. 
If however $\ell$ is odd, then~\eqref{eq:gamma_m2} and~\eqref{eq:links} imply that
\begin{equation}
\label{eq:jumps2}
	\vv{m-1,2(m-1) + 4j_0 + \delta} \sim \vv{0,\ell - bj_0-\delta a}
\end{equation}
holds for each $j_0$ with $0 \leq j_0 < n_0$ and each $\delta \in \{0,1\}$.

\subsection{The construction}

We have thus seen that (at least in the case of $n > 8$ or $n = 8$ with $(a,b) \in \{(1,4), (5,4)\}$) the graph $\G$ (and the action of $\gamma$) is completely determined by the parameters $m$, $n$, $a$, $b$ and $\ell$. We now record the construction obtained in this way for future reference.

\begin{construction}
\label{con:Xb}
Let $m \geq 3$ and $n \geq 8$ be integers where $n = 4n_0$ for some integer $n_0 \geq 2$. Furthermore, let $a,b,\ell \in \ZZ_n$ be such that~\eqref{eq:b} and~\eqref{eq:b_a_cond} hold and that one of $\ell - 2$ and $\ell - 3$ is divisible by $4$. Then the graph $\cXb(m,n,a,b,\ell)$ is the graph with vertex set $\{\vv{i,j} \colon i \in \ZZ_m,\ j \in \ZZ_n\}$ such that each $\vv{i,j}$ is adjacent to $\vv{i,j\pm 1}$ and we also have the adjacencies from~\eqref{eq:links}, together with all the adjacencies from~\eqref{eq:jumps1} or~\eqref{eq:jumps2}, depending on whether $\ell$ is even or odd, respectively.
\end{construction}

The example $\cXb(5,12,1,8,7)$ is given in Figure~\ref{fig:example_(5,12,1,8,7)}. It turns out that this graph has a regular automorphism group, and that the ``natural'' $2$-factor consisting of the five $12$-cycles is preserved by the full automorphism group. 
%
\begin{figure}[h]
\begin{center}
\subfigure
{
\begin{tikzpicture}[scale = .4]
\foreach \i in {0,1,2,3,4}{
\foreach \j in {0,1,2,3,4,5,6,7,8,9,10,11}{
\node[vtx, fill = white, inner sep = 1.5pt,] (A\i\j) at (3*\i, 1.6*\j) {};
}}
\begin{scope}[on background layer]
\foreach \j in {0,1,2,3,4,5,6,7,8,9,10}{
\foreach \i in {0,1,2,3,4}{
\draw[thick] let \n1 = {int(mod(\j+1, 12))} in (A\i\j) -- (A\i\n1);
}}
\foreach \i in {0,1,2,3,4}{
\draw[thick] (A\i 0) to[bend left = 10] (A\i 11);
}
\foreach \i in {0,1,2,3}{
\draw[thick] let \n1 = {int(mod(\i+1, 5))}, \n2 = {int(mod(2*\i, 12))} in (A\i\n2) -- (A\n1\n2);
\draw[thick] let \n1 = {int(mod(\i+1, 5))}, \n2 = {int(mod(2*\i + 1, 12))} in (A\i\n2) -- (A\n1\n2);
\draw[thick] let \n1 = {int(mod(\i+1, 5))}, \n2 = {int(mod(2*\i + 4, 12))}, \n3 = {int(mod(2*\i + 8, 12))} in (A\i\n2) -- (A\n1\n3);
\draw[thick] let \n1 = {int(mod(\i+1, 5))}, \n2 = {int(mod(2*\i + 5, 12))}, \n3 = {int(mod(2*\i + 9, 12))} in (A\i\n2) -- (A\n1\n3);
\draw[thick] let \n1 = {int(mod(\i+1, 5))}, \n2 = {int(mod(2*\i + 8, 12))}, \n3 = {int(mod(2*\i + 4, 12))} in (A\i\n2) -- (A\n1\n3);
\draw[thick] let \n1 = {int(mod(\i+1, 5))}, \n2 = {int(mod(2*\i + 9, 12))}, \n3 = {int(mod(2*\i + 5, 12))} in (A\i\n2) -- (A\n1\n3);
}
\draw[thick] (A02) -- (A45);
\draw[thick] (A03) -- (A44);
\draw[thick] (A06) -- (A49);
\draw[thick] (A07) -- (A48);
\draw[thick] (A010) -- (A41);
\draw[thick] (A011) -- (A40);
\end{scope}
\end{tikzpicture}
}
\quad \quad
\subfigure
{
\begin{tikzpicture}[scale = .4]
\foreach \i in {0,1,2,3,4}{
\pgfmathtruncatemacro{\rad}{\i+1};
\foreach \j in {0,1,2,3,4,5,6,7,8,9,10,11}{
\node[vtx, fill = white, inner sep = 1.5pt,] (A\i\j) at (360*\j/12:1.7*\rad) {};
}}
\begin{scope}[on background layer]
\foreach \j in {0,1,2,3,4,5,6,7,8,9,10,11}{
\foreach \i in {0,1,2,3,4}{
\draw[thick] let \n1 = {int(mod(\j+1, 12))} in (A\i\j) -- (A\i\n1);
}}
\foreach \i in {0,1,2,3}{
\draw[thick] let \n1 = {int(mod(\i+1, 5))}, \n2 = {int(mod(2*\i, 12))} in (A\i\n2) -- (A\n1\n2);
\draw[thick] let \n1 = {int(mod(\i+1, 5))}, \n2 = {int(mod(2*\i + 1, 12))} in (A\i\n2) -- (A\n1\n2);
\draw[thick] let \n1 = {int(mod(\i+1, 5))}, \n2 = {int(mod(2*\i + 4, 12))}, \n3 = {int(mod(2*\i + 8, 12))} in (A\i\n2) -- (A\n1\n3);
\draw[thick] let \n1 = {int(mod(\i+1, 5))}, \n2 = {int(mod(2*\i + 5, 12))}, \n3 = {int(mod(2*\i + 9, 12))} in (A\i\n2) -- (A\n1\n3);
\draw[thick] let \n1 = {int(mod(\i+1, 5))}, \n2 = {int(mod(2*\i + 8, 12))}, \n3 = {int(mod(2*\i + 4, 12))} in (A\i\n2) -- (A\n1\n3);
\draw[thick] let \n1 = {int(mod(\i+1, 5))}, \n2 = {int(mod(2*\i + 9, 12))}, \n3 = {int(mod(2*\i + 5, 12))} in (A\i\n2) -- (A\n1\n3);
}
\draw[thick] (A02) -- (A45);
\draw[thick] (A03) -- (A44);
\draw[thick] (A06) -- (A49);
\draw[thick] (A07) -- (A48);
\draw[thick] (A010) -- (A41);
\draw[thick] (A011) -- (A40);
\end{scope}
\end{tikzpicture}
}
\caption{Two presentations of the graph $\cXb(5,12,1,8,7)$.}
\label{fig:example_(5,12,1,8,7)}
\end{center}
\end{figure}
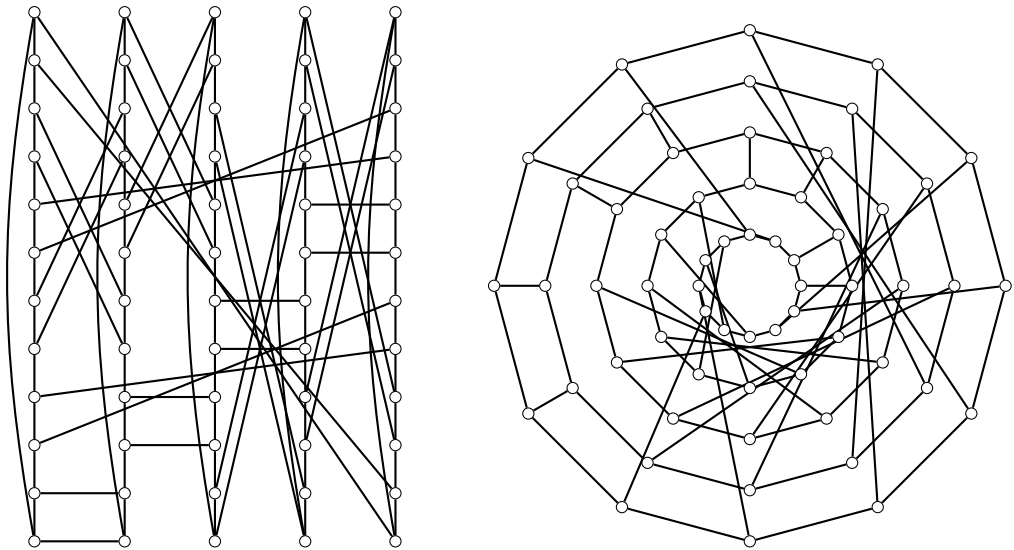

We can now also state the main result of this section. At the same time we make the agreement that whenever we are working with a graph $\G = \cXb(m,n,a,b,\ell)$ from Construction~\ref{con:Xb}, we let $V_i = \{\vv{i,j} \colon j \in \ZZ_n\}$ for all $i \in \ZZ_m$, and we let $C_i$, $i \in \ZZ_m$, be the subgraphs of $\G$ induced on $V_i$ (which are thus cycles of length $n$). We let $\cC = \{C_i \colon i \in \ZZ_m\}$, to which we refer as {\em the natural} {\em $2$-factor} of $\G$. We also call the edges arising from~\eqref{eq:links} the {\em links} of $\G$ and those arising from~\eqref{eq:jumps1} or~\eqref{eq:jumps2} (depending on the parity of $\ell$) the {\em non-links} of $\G$.

\begin{proposition}
\label{pro:the_Xb_graphs}
Let $\G$ be a cubic graph admitting a vertex-transitive subgroup $G$ of the automorphism group $\Aut(\G)$ and a $G$-invariant partition of the edge set of $\G$ into a $2$-factor $\cC$ and a $1$-factor such that $\G$ is of bialternating cycle quotient type with respect to $\cC$. Let $m = |\cC|$ and let $n$ be the length of the cycles in $\cC$. Then $n$ is divisible by $4$. Moreover, if $n > 8$ or $n = 8$ and the subgraph induced on the vertices of two adjacent $8$-cycles from $\cC$ is bipartite, then $\G$ is isomorphic to some $\cXb(m,n,a,b,\ell)$ from Construction~\ref{con:Xb}, where the cycles of $\cC$ correspond to the subgraphs induced on the sets $V_i = \{\vv{i,j} \colon j \in \ZZ_n\}$, $i \in \ZZ_m$.
\end{proposition}

\section{Existence of a vertex-transitive subgroup preserving the partition}
\label{sec:VT}

Proposition~\ref{pro:the_Xb_graphs} implies that, except for the cases when the cycles from the $2$-factor $\cC$ are $8$-cycles and the subgraph induced on the vertices of two adjacent $8$-cycles from $\cC$ is not bipartite, it suffices to analyze the graphs $\G = \cXb(m,n,a,b,\ell)$ from Construction~\ref{con:Xb}. It should not come as a surprise that not all graphs from Construction~\ref{con:Xb} admit a vertex-transitive subgroup $G \leq \Aut(\G)$ which preserves the natural $2$-factor $\cC$ of $\G$. In fact, some graphs from Construction~\ref{con:Xb} are not even vertex-transitive. For instance, using suitable software one can verify that the graphs $\cXb(3,12,1,4,2)$ and $\cXb(3,12,1,4,3)$ have automorphism groups of orders $18$ and $12$, respectively. On the other hand, $\cXb(3,12,1,4,10)$ is vertex-transitive and has an automorphism group of order $72$. 

The aim of this section is to obtain a necessary and sufficient condition for a graph $\G = \cXb(m,n,a,b,\ell)$ from Construction~\ref{con:Xb} to admit a vertex-transitive subgroup $G \leq \Aut(\G)$ such that the natural $2$-factor $\cC$ of $\G$ is $G$-invariant. Throughout the section let $\G = \cXb(m,n,a,b,\ell)$ be a graph from Construction~\ref{con:Xb} and recall that we assume by definition that~\eqref{eq:b} and~\eqref{eq:b_a_cond} both hold and that one of $\ell - 2$ and $\ell - 3$ is divisible by $4$. 

\subsection{Existence of the automorphism $\gamma$}

Let $\gamma$ be the permutation of the vertex set of $\G$, defined by~\eqref{eq:gamma} and~\eqref{eq:gamma_m1} when $\ell$ is even, and by~\eqref{eq:gamma} and~\eqref{eq:gamma_m2} when $\ell$ is odd. We now determine a necessary and sufficient condition for $\gamma$ to be an automorphism of $\G$. It is clear that $\gamma$ preserves the set of non-links of $\G$, and therefore also preserves the natural $2$-factor $\cC$. In view of~\eqref{eq:gamma} and~\eqref{eq:links}, it is also clear that for each $i$, $0 \leq i \leq m-3$, it maps the links between the sets $V_i$ and $V_{i+1}$ to the links between the sets $V_{i+1}$ and $V_{i+2}$. We thus only need to consider the links between $V_{m-2}$ and $V_{m-1}$, and the links between $V_{m-1}$ and $V_0$. 

Let $\varepsilon$ be $1$ or $-1$, depending on whether $\ell$ is even or odd, respectively. Pick a $\delta \in \{0,1\}$ and an integer $j_0$ with $0 \leq j_0 < n_0$. By~\eqref{eq:gamma},~\eqref{eq:gamma_m1} and~\eqref{eq:gamma_m2}, the vertex $\vv{m-2,2(m-2)+4j_0+\delta}$ is mapped to $\vv{m-1,2(m-1)+4j_0+\delta}$, while its outside neighbor $\vv{m-1,2(m-2)+bj_0+\delta a}$ is mapped to $\vv{0,\ell + \varepsilon (bj_0+\delta a)}$. Therefore,~\eqref{eq:jumps1} and~\eqref{eq:jumps2} show that $\gamma$ maps the links between $V_{m-2}$ and $V_{m-1}$ to links of $\G$. 
Finally, $\gamma$ maps the adjacent vertices $\vv{m-1,2(m-1)+4j_0+\delta}$ and $\vv{0,\ell+\varepsilon(bj_0+\delta a)}$ to 
\[
	v = \vv{0,\ell+\varepsilon(4j_0+\delta+2)} \quad \text{and} \quad w = \vv{1, \ell+2+\varepsilon(bj_0+\delta a)},
\]
respectively. 

Suppose first that $\varepsilon = 1$ (in which case $\ell = 4\ell_0+2$ is even). Let us write the second component of the vertex $v$ as $4(\ell_0 + 1 + j_0) + \delta$. Then~\eqref{eq:links} implies that $\No{v} = \vv{1,b(\ell_0 + 1 + j_0) + \delta a}$, and so $v$ and $w$ are adjacent if and only if 
$$
	\ell + 2 + bj_0 + \delta a = b(\ell_0 + 1 + j_0) + \delta a = b_0(4\ell_0 + 4) + bj_0 + \delta a
$$
holds in $\ZZ_n$, which is equivalent to the condition
\begin{equation}
\label{eq:cond_gamma1}
	(\ell + 2)(b_0 - 1) = 0\ \text{in}\ \ZZ_n.
\end{equation}

Suppose now that $\varepsilon = -1$ (in which case $\ell = 4\ell_0+3$ is odd). We now write the second component of $v$ as $4(\ell_0-j_0) + 1-\delta$, and find that $\No{v} = \vv{1,b(\ell_0-j_0)+(1-\delta)a}$. The vertices $v$ and $w$ are thus adjacent if and only if 
$$
	\ell + 2 - bj_0 - \delta a = b(\ell_0-j_0) + (1-\delta)a
$$
holds in $\ZZ_n$, which is equivalent to the condition
\begin{equation}
\label{eq:cond_gamma2}
	\ell + 2 = b\ell_0 + a\ \text{in}\ \ZZ_n.
\end{equation}
Suppose~\eqref{eq:cond_gamma2} does indeed hold and note that in this case~\eqref{eq:b_a_cond} implies that in $\ZZ_n$
$$
	\ell - 3 = 4\ell_0 = b_0(b\ell_0) = b_0(\ell+2-a) = 4b_0(\ell_0+1-a_0) = b\ell_0 + b - ba_0 = \ell + 2 - a + b + a - 1.
$$ 
Therefore, $b+4 = 0$ must hold. Plugging this into~\eqref{eq:cond_gamma2} we then also get $a = 2\ell - 1$. Of course, if $b = -4$ and $a = 2\ell-1$, then~\eqref{eq:cond_gamma2} clearly holds. This finally proves the following result.

\begin{lemma}
\label{le:gamma}
Let $\G = \cXb(m,n,a,b,\ell)$ be a graph from Construction~\ref{con:Xb}. Then there exists an automorphism $\gamma$ of $\G$ preserving the natural $2$-factor $\cC$ and mapping the vertex $\vv{0,0}$ to the vertex $\vv{1,2}$, if and only if one of the following holds:
\begin{itemize}
\item $\ell$ is odd and $b = -4$ and $a = 2\ell-1$ both hold in $\ZZ_n$, or
\item $\ell$ is even and $(\ell + 2)(b_0 - 1) = 0$ holds in $\ZZ_n$, where $b = 4b_0$.
\end{itemize}
\end{lemma}

\subsection{Existence of the automorphism $\alpha$}

We next determine a necessary and sufficient condition for a graph $\G = \cXb(m,n,a,b,\ell)$, satisfying the conditions from Lemma~\ref{le:gamma}, to admit an automorphism $\alpha$ preserving the natural $2$-factor $\cC$ and such that its action on the sets $V_0$ and $V_1$ is as in~\eqref{eq:alpha}. As we now demonstrate, the requirement that $\alpha$ preserves $\cC$ and maps the vertices of $V_0 \cup V_1$ as in~\eqref{eq:alpha}, completely determines the action of $\alpha$. We thus only need to determine this action and then establish a necessary and sufficient condition under which the obtained permutation is indeed an automorphism of $\G$. To start, observe that using the idea of the proof of Lemma~\ref{le:G_regular} enables us to show that the only automorphism of $\G$, fixing a vertex and preserving the natural $2$-factor $\cC$, is the identity. Therefore, $\alpha$ (if it exists) is an involution which interchanges each $V_i$ with $V_{1-i}$, where the indices are computed modulo $m$. 

Determining the exact action of $\alpha$ is somewhat tedious. We consider the two cases depending on the parity of $\ell$ separately. We provide most of the details for the first case and leave some for the second one to the reader.
\medskip

\noindent {\sc Case 1:} $\ell = 4\ell_0 + 2$ for some $0 \leq \ell_0 < n_0$.\\
Lemma~\ref{le:gamma} implies that in this case $4(\ell_0 + 1)(b_0 - 1) = 0$ holds in $\ZZ_n$, which by~\eqref{eq:b} is equivalent to $b\ell_0 = \ell + 2 - b$. For future reference we record this fact and an immediate consequence of~\eqref{eq:b_a_cond}:
\begin{equation}
\label{eq:alpha_aux}
	b\ell_0 = \ell + 2 - b\ \text{and}\ \delta ba_0 = \delta - \delta a\ \text{for\ each}\ \delta \in \{0,1\}.
\end{equation}
Let $0 \leq j_0 < n_0$ and $\delta \in \{0,1\}$. By~\eqref{eq:jumps1} we must have that $\alpha(\vv{m-1, 2(m-1) + 4j_0 + \delta}) = \No{(\alpha(\vv{0, \ell + bj_0 + \delta a}))}$. Combining together~\eqref{eq:alpha},~\eqref{eq:b_a_cond},~\eqref{eq:links} and~\eqref{eq:alpha_aux} we obtain 
$$
	\No{(\alpha(\vv{0, \ell + bj_0 + \delta a}))} = \No{\vv{1, \ell + bj_0 + \delta a}} = \No{\vv{1, 2 + 4(\ell_0 + b_0j_0 + \delta a_0) + \delta}} =  \vv{2,\ell+4-b+4j_0+\delta},
$$
and so $\alpha$ maps the links connecting the vertices of $V_0$ and $V_{m-1}$ to edges of $\G$ if and only if 
\begin{equation}
\label{eq:alpha_aux1}
	\alpha(\vv{m-1, 2(m-1)+j}) = \vv{2, \ell + 4 - b + j}\ \text{for\ all}\ j \in \ZZ_n.
\end{equation}
Next, $\vv{m-2, 2(m-2) + 4j_0 + \delta} \sim \vv{m-1, 2(m-1) - 2 + bj_0 + \delta a}$, the latter of which is mapped by $\alpha$ to $\vv{2, \ell + 2 - b + bj_0 + \delta a}$. By~\eqref{eq:b_a_cond},~\eqref{eq:links} and~\eqref{eq:alpha_aux} we find that
$$
	\No{\vv{2, \ell + 2 - b + bj_0 + \delta a}} = \No{\vv{2, 4 + 4(\ell_0 - b_0 + b_0j_0 + \delta a_0) + \delta}} = \vv{3, 4 + b(\ell_0 - b_0 + b_0j_0 + \delta a_0) + \delta a} = \vv{3, \ell+2-b + 4j_0 + \delta}.
$$
Thus, $\alpha$ maps the links connecting the vertices of $V_{m-2}$ and $V_{m-1}$ to edges of $\G$ if and only if 
$$
	\alpha(\vv{m-2, 2(m-2) + j}) = \vv{3, \ell + 2 - b + j}\ \text{for\ all}\ j \in \ZZ_n.
$$
Continuing in this way we can determine the action of the permutation $\alpha$ completely. We find that for all integers $i$ and $\varepsilon \in \{0,1\}$ such that $1 \leq 2i-\varepsilon \leq m-2$ and for all $j \in \ZZ_n$, the following holds:
\begin{equation}
\label{eq:alpha_fin}
	\alpha(\vv{m-2i+\varepsilon, 2(m-2i+\varepsilon)+j}) = \vv{2i-\varepsilon+1, \ell + 2 + 2\varepsilon - i b + j}.
\end{equation}
Moreover, a permutation $\alpha$ defined by~\eqref{eq:alpha} and~\eqref{eq:alpha_fin} maps the non-links of $\G$ to edges (in fact to non-links) of $\G$ and also maps all links of $\G$, except possibly those connecting the vertices of $V_1$ and $V_2$, to edges of $\G$. It thus remains to consider the links connecting the vertices of $V_1$ and $V_2$. To be able to do this we need to distinguish two subcases depending on the parity of $m$.
\medskip 

\noindent {\sc Subcase 1.1:} $m$ is odd, that is $m = 2m_0 + 1$ for some $m_0 \geq 1$.\\
Recall that $\No{\vv{1,2+4j_0+\delta}} = \vv{2,2+bj_0+\delta a}$ for each $j_0$ and $\delta$ with $0 \leq j_0 < n_0$ and $\delta \in \{0,1\}$. Since $2 = m - 2m_0 + 1$,~\eqref{eq:b_a_cond} and \eqref{eq:alpha_fin} imply that 
$$
	\alpha(\vv{2,2+bj_0+\delta a}) = \vv{m-1,\ell+4-m_0 b + bj_0 + \delta a - 2} = \vv{m-1, 2(m-1) + 4(-m_0 + \ell_0 + 1 - m_0 b_0 + b_0j_0 + \delta a_0) + \delta}.
$$
According to~\eqref{eq:jumps1}, where we also use~\eqref{eq:b_a_cond} and~\eqref{eq:alpha_aux}, the outside neighbor of this vertex is
$$
	\vv{0,\ell - bm_0 + \ell + 2 - b + b - 4m_0 + 4j_0 + \delta} = \vv{0, 2 + 4j_0 + \delta + 2\ell - bm_0 - 4m_0}.
$$
Since $\alpha(\vv{1,2+4j_0+\delta}) = \vv{0,2+4j_0+\delta}$, this finally implies that $\alpha$ maps the links between the vertices of $V_1$ and $V_2$ to edges of $\G$ if and only if 
\begin{equation}
\label{eq:l_even_m_odd}
	2\ell = m_0(b+4)\ \text{in}\ \ZZ_n.
\end{equation}
By~\eqref{eq:b_a_cond}, the right-hand side of~\eqref{eq:l_even_m_odd} does not change when we multiply it by $b_0$, and so 
$$ 
	2\ell = 2\ell b_0 = 2(4\ell_0+2)b_0 = 2b\ell_0+b = 2(\ell+2-b)+b = 2\ell + 4 - b.
$$
Therefore, $b = 4$. Plugging this into~\eqref{eq:l_even_m_odd}, we find that $8m_0 = 2\ell = 8\ell_0 + 4$, and so $n \equiv 4 \pmod{8}$ and $\ell = 2n_0 + 4m_0 = 2n_0 + 2m-2$. Finally, by~\eqref{eq:b_a_cond} we have that $4a_0 + a = 1$, that is $2(a-1) = 0$. Since $a-1$ is divisible by $4$, this thus implies that $a = 1$. To summarize, if $\ell$ is even and $m$ is odd, then the automorphism $\alpha$ exists if and only if $n \equiv 4 \pmod{8}$, $a = 1$, $b = 4$ and $\ell = n/2 + 2m - 2$. 
\medskip

\noindent {\sc Subcase 1.2:} $m = 2m_0$ is even.\\
This time $2 = m - 2(m_0-1)$, and so~\eqref{eq:alpha_fin} implies that
$$
	\alpha(\vv{2,2+bj_0+\delta a}) = \vv{m-1,\ell+2-(m_0-1) b + bj_0 + \delta a - 2} = \vv{m-1, 2(m-1) + 4(-m_0 + \ell_0 + 1 - (m_0-1) b_0 + b_0j_0 + \delta a_0) + \delta}.
$$
Similarly as before we find that the outside neighbor of this vertex is $\vv{0,2\ell-bm_0-4m_0+4+2+4j_0+\delta}$, implying that $\alpha$ maps the links between the vertices of $V_1$ and $V_2$ to edges of $\G$ if and only if 
\begin{equation}
\label{eq:l_even_m_even}
	2(\ell + 2) = m_0(b+4)\ \text{in}\ \ZZ_n.
\end{equation}
\medskip

\noindent {\sc Case 2:} $\ell = 4\ell_0 + 3$ for some $0 \leq \ell_0 < n_0$.\\
Recall that Lemma~\ref{le:gamma} forces $b = -4$ and $a = 2\ell-1$. We first notice that this time 
\begin{equation}
\label{eq:alpha_aux2}
	-4(\ell_0-\delta a_0) + (1-\delta)a = 3- \ell + \delta (a-1) + a - \delta a = 2 + \ell - \delta 
\end{equation}
holds for each $\delta \in \{0,1\}$. We now proceed similarly as in Case 1. 
Let $0 \leq j_0 < n_0$ and $\delta \in \{0,1\}$. This time~\eqref{eq:jumps2} gives $\vv{m-1, 2(m-1) + 4j_0 +  \delta} \sim \vv{0, \ell + 4j_0 - \delta a}$. The latter of these two vertices is mapped by $\alpha$ to $\vv{1, \ell + 4j_0 - \delta a}$. Moreover,~\eqref{eq:alpha_aux2} yields
$$
	\No{\vv{1, \ell + 4j_0 - \delta a}} = \No{\vv{1, 2 + 4(\ell_0 + j_0 - \delta a_0) + 1 - \delta}} = \vv{2, 2 - 4(\ell_0 - \delta a_0) + (1-\delta)a - 4j_0} = \vv{2, 4 + \ell - 4j_0 - \delta}. 
$$
Therefore, $\alpha$ maps the links connecting the vertices of $V_0$ and $V_{m-1}$ to edges of $\G$ if and only if 
$$
	\alpha(\vv{m-1, 2(m-1) + j}) = \vv{2, \ell + 4 - j}\ \text{for\ all}\ j \in \ZZ_n.
$$
Just like in Case~1, continuing in this way we find that the permutation $\alpha$ must act in such a way that for all integers $i$ and $\varepsilon \in \{0,1\}$ such that $1 \leq 2i-\varepsilon \leq m-2$ and for all $j \in \ZZ_n$, the following holds:
\begin{equation}
\label{eq:alpha_fin2}
	\alpha(\vv{m-2i+\varepsilon, 2(m-2i+\varepsilon)+j}) = \vv{2i-\varepsilon+1, \ell + 4i + 2(\varepsilon-1) - j}.
\end{equation}
Moreover, the permutation $\alpha$, defined by~\eqref{eq:alpha} and~\eqref{eq:alpha_fin2}, maps the edges of $\G$ to edges of $\G$, except possibly for the links connecting the vertices of $V_1$ and $V_2$. Unlike in Case~1, where both parities of $m$ were possible, this time $m$ has to be odd. To see this, suppose to the contrary that $m = 2m_0$ and note that $2 = m -2(m_0-1)$. Then~\eqref{eq:alpha_fin2} implies that
$$
	\alpha^2(\vv{m-1,2(m-1)}) = \alpha(\vv{2,\ell + 4}) = \vv{m-1, \ell + 4(m_0-1) - 2 - \ell} = \vv{m-1, 2(m-1)-4}, 
$$
contradicting the fact that $n > 4$ and that $\alpha$ is an involution. Therefore, $m = 2m_0 + 1$, as claimed. We leave the verification that in this case the permutation $\alpha$, defined by~\eqref{eq:alpha} and~\eqref{eq:alpha_fin2}, does indeed map the links connecting the vertices of $V_1$ and $V_2$ to edges of $\G$ to the reader. This finally proves the following result.

\begin{lemma}
\label{le:alpha}
Let $\G = \cXb(m,n,a,b,\ell)$ be a graph from Construction~\ref{con:Xb}, where the parameters also satisfy the conditions of Lemma~\ref{le:gamma}. Then there exists an automorphism $\alpha$ of $\G$ interchanging each $V_i$ with $V_{1-i}$ (indices being computed modulo $m$) and mapping the vertex $\vv{0,0}$ to $\vv{1,0}$, if and only if one of the following holds:
\begin{itemize}\itemsep = 0pt
	\item $m$ and $\ell$ are both odd, and $b = -4$ and $a = 2\ell-1$ both hold in $\ZZ_n$;	
	\item $m$ is odd, $\ell$ is even, $a = 1$, $b = 4$, $n = 4n_0$ for an odd integer $n_0 \geq 3$, and $\ell = 2n_0 + 2m - 2$ holds in $\ZZ_n$.
	\item $m$ and $\ell$ are both even and \eqref{eq:l_even_m_even} holds;
\end{itemize}
\end{lemma}

\subsection{Sufficiency for vertex transitivity}

Lemmas~\ref{le:gamma} and~\ref{le:alpha} give a necessary and sufficient condition for a graph $\G$ from Construction~\ref{con:Xb} to admit (the unique) automorphisms $\alpha$ and $\gamma$, preserving the natural $2$-factor $\cC$ and such that $\alpha$ and $\gamma$ map $\vv{0,0}$ to $\vv{1,0}$ and $\vv{1,2}$, respectively. The aim of this subsection is to show that this condition is in fact sufficient for the existence of a regular subgroup of $\Aut(\G)$ preserving $\cC$.  

To this end let $\G = \cXb(m,n,a,b,\ell)$ be a graph satisfying all the assumptions and conditions of Lemma~\ref{le:alpha} and let $\alpha$ and $\gamma$ be as in the previous two subsections. Using~\eqref{eq:alpha} and~\eqref{eq:gamma} we obtain
$$
	(\alpha\gamma)(\vv{0,0}) = \alpha(\vv{1,2}) = \vv{0,2}. 
$$
Since the group $\langle \alpha, \gamma \rangle$ preserves the natural $2$-factor $\cC$, the automorphism $\alpha\gamma$ thus acts as a $2$-step rotation on $C_0$ (it acts according to~\eqref{eq:rho_V_0}). The fact that $\gamma$ cyclically permutes the $m$ cycles from $\cC$ therefore implies that the group $\la \alpha, \gamma \ra$ has at most two orbits on the vertex set of $\G$. In fact, the nature of the action of $\gamma$ reveals that the $\la \alpha, \gamma \ra$-orbit of $\vv{0,0}$ for sure contains all the vertices $\vv{i,j}$ with $j$ even. By~\eqref{eq:gamma_m2} we thus see that if $\ell$ is odd, then $\G$ is actually $\la \alpha, \gamma \ra$-vertex-transitive (and consequently $\la \alpha, \gamma \ra$ is regular).

It thus suffices to consider the case that $\ell = 4\ell_0 + 2$ is even. If $m$ is odd, then Lemma~\ref{le:alpha} implies that $a = 1$, $b = 4$ and $\ell = 2n_0 + 2m - 2$, where $n = 4n_0$ for some odd integer $n_0 \geq 3$. Then~\eqref{eq:jumps1} becomes
$$
	\vv{m-1,2(m-1)+4j_0+\delta} \sim \vv{0, 2(m-1) + 2n_0 + 4j_0 + \delta},
$$ 
and so~\eqref{eq:links} clearly implies that the permutation $\beta$, mapping each $\vv{i,j}$ to $\vv{i,1-j}$, is an automorphism of $\G$. Of course, $\la \alpha, \beta, \gamma \ra$ is then a vertex-transitive subgroup of $\Aut(\G)$ and it preserves the $2$-factor $\cC$.

This leaves us with the possibility that $m = 2m_0$ for some $m_0 \geq 2$ (recall that we are also assuming $\ell = 4\ell_0 + 2$). Recall that Lemmas~\ref{le:gamma} and~\ref{le:alpha} imply that
\begin{equation}
\label{eq:beta_aux}
	(\ell+2)(b_0-1) = 0\ \text{and}\ 2(\ell + 2) = m_0(b+4) = 2m(b_0+1)\ \text{in}\ \ZZ_n.
\end{equation}
Let $\beta$ be the permutation of the vertex set of $\G$, defined by setting 
\begin{equation}
\label{eq:beta}
	\beta(\vv{2i + \varepsilon, j}) = \vv{2i + \varepsilon, \varepsilon (a-1) + 1+ i(4-b)-j}\ \text{for\ each}\ i,\ 0 \leq i < m_0,\ \varepsilon \in \{0,1\}\ \text{and}\ j \in \ZZ_n.
\end{equation}
It is clear that $\beta$ is a permutation of $V(\G)$ and that it preserves the set of the non-links of $\G$. We now show that it also preserves the set of the links of $\G$. For each $i$ with $0 \leq i < m_0$, each $j_0$ with $0 \leq j_0 < n_0$ (recall that $n = 4n_0$) and each $\delta \in \{0,1\}$ consider the adjacency
$$
	\vv{2i, 4i+4j_0+\delta} \sim \vv{2i+1,4i+bj_0+\delta a}.
$$
By~\eqref{eq:b_a_cond} and~\eqref{eq:beta} we first establish that
$$
	\beta(\vv{2i, 4i+4j_0+\delta}) = \vv{2i, 1+i(4-b)-4i-4j_0-\delta} = \vv{2i, 4i+4(-b_0 i - i - j_0) + 1-\delta},
$$
which by~\eqref{eq:b_a_cond} and~\eqref{eq:links} is adjacent to the vertex $\vv{2i+1,-bi-bj_0+(1-\delta)a}$. Since 
$$
	\beta(\vv{2i+1,4i+bj_0+\delta a}) = \vv{2i+1,a+i(4-b)-4i-bj_0-\delta a},
$$
this thus shows that $\beta$ preserves the set of links connecting the vertices of $V_{2i}$ and $V_{2i+1}$ for each $i$ with $0 \leq i < m_0$. Similarly, consider the adjacency 
$$
	\vv{2i+1,4i+2+4j_0+\delta} \sim \vv{2i+2, 4i+2+bj_0+\delta a},
$$
where $0 \leq i \leq m_0-2$, $0 \leq j_0 < n_0$ and $\delta \in \{0,1\}$. By~\eqref{eq:beta}, the first of these vertices is mapped by $\beta$ to the vertex
$$
	\vv{2i+1,a+i(4-b)-4i-2-4j_0-\delta} = \vv{2i+1,4i+2+4(a_0-1-b_0 i-i-j_0)+1-\delta}.
$$
By~\eqref{eq:b_a_cond} and~\eqref{eq:links}, its outside neighbor is
$$
	\vv{2i+2,4i+2+ba_0-b-4i-bi-bj_0+(1-\delta)a} = \vv{2i+2,3-b(1+i+j_0)-\delta a}.
$$  
Since~\eqref{eq:beta} also implies that 
$$
	\beta(\vv{2i+2, 4i+2+bj_0+\delta a}) = \vv{2i+2, 1 + (i+1)(4-b)-4i-2-bj_0-\delta a} = \vv{2i+2, 3 - b - bi - bj_0 - \delta a},
$$
this shows that $\beta$ also preserves the set of the links connecting the vertices of $V_{2i+1}$ and $V_{2i+2}$ for each $i$ with $0 \leq i \leq m_0-2$. Therefore, $\beta$ is an automorphism of $\G$ if and only if it preserves the set of the links  
$$
	\vv{m-1, 2(m-1)+4j_0+\delta} \sim \vv{0, \ell+bj_0+\delta a},
$$ 
connecting the vertices of $V_{m-1}$ to those of $V_0$. The first of these two vertices is mapped by $\beta$ to
$$
	\vv{m-1,a+(m_0-1)(4-b)-2(m-1)-4j_0-\delta} = \vv{m-1, 2(m-1) + 4(-m_0+a_0-b_0m_0+b_0-j_0)+1-\delta}.
$$
By~\eqref{eq:b_a_cond},~\eqref{eq:jumps1} and~\eqref{eq:beta_aux}, the outside neighbor of this vertex is  
$$
	\vv{0,\ell-bm_0+1-a-4m_0+4-bj_0+(1-\delta)a} = \vv{0, \ell - 2(\ell + 2) + 5 - bj_0 - \delta a} = \vv{0,1-\ell - bj_0 - \delta a}.
$$
Since this is precisely the vertex $\beta(\vv{0, \ell+bj_0+\delta a})$, this finally proves that $\beta$ is an automorphism of $\G$. We can now state the main result of this section.

\begin{proposition}
\label{pro:VT2}
Let $\G$ be a cubic vertex-transitive graph admitting a partition of its edge set into a $2$-factor $\cC$ and a $1$-factor such $\G$ is of bialternating cycle quotient type with respect to $\cC$. Let $n$ be the length of the cycles from $\cC$ and suppose that either $n > 8$, or $n = 8$ and the subgraph of $\G$ induced on two adjacent $8$-cycles from $\cC$ is bipartite. Then $\G$ admits a vertex-transitive subgroup $G$ of $\Aut(\G)$ preserving the $2$-factor $\cC$, if and only if $n$ is divisible by $4$ and $\G$ is isomorphic to a graph $\cXb(m,n,a,b,\ell)$ from Construction~\ref{con:Xb}, where the cycles of $\cC$ correspond to the subgraphs of $\G$ induced on the sets $V_i = \{\vv{i,j} \colon j \in \ZZ_n\}$, $i \in \ZZ_m$, and one of the following holds:
\begin{itemize}\itemsep = 0pt
	\item $m$ and $\ell$ are both odd, and $b = -4$ and $a = 2\ell-1$ both hold in $\ZZ_n$;
	\item $m$ is odd, $\ell$ is even, $a = 1$, $b = 4$, $n = 4n_0$ for an odd integer $n_0 \geq 3$, and $\ell = 2n_0 + 2m - 2$ holds in $\ZZ_n$.
	\item $m$ and $\ell$ are both even and $(\ell + 2)(b_0 - 1) = 0$ and $2(\ell + 2) = 2m(b_0+1)$ both hold in $\ZZ_n$, where $b = 4b_0$.
\end{itemize}
Moreover, $G$ is a regular group and $\G$ is a Cayley graph of $G$ with respect to a connection set consisting of three involutions.
\end{proposition}

\begin{proof}
The only claim that has not been proved thus far is that if a vertex-transitive group $G$ preserves the $2$-factor $\cC$, then $\G$ is a Cayley graph of $G$ with respect to a connection set consisting of three involutions. Since $G$ is regular by Lemma~\ref{le:G_regular}, we only need to see that each of the three elements of $G$, mapping $\vv{0,0}$ to one of its neighbors, maps this neighbor back to $\vv{0,0}$. This follows from the fact that $\No{\vv{0,0}}, \No{\vv{0,1}} \in V_1$, while $\No{\vv{0,-1}}, \No{\vv{0,-2}}, \No{\vv{0,2}} \notin V_1$. In particular, the element of $G$, mapping $\vv{0,0}$ to $\vv{0,1}$, preserves $V_0$ and $V_1$, and thus cannot map $\vv{0,1}$ to $\vv{0,2}$. Similarly, the element of $G$, mapping $\vv{0,0}$ to $\vv{0,-1}$, preserves $V_0$ but interchanges $V_1$ with $V_{m-1}$, and thus cannot map $\vv{0,-1}$ to $\vv{0,-2}$.
\end{proof}

\section{The examples with $n = 8$}
\label{sec:n=8}

In Subsection~\ref{subsec:n=8}, we pointed out that when the $2$-factor $\cC$ consists of $8$-cycles, the situation is a bit more involved. The aim of this section is to analyze this possibility. 

Throughout this section we therefore assume that $\G$ is a cubic graph admitting a vertex-transitive subgroup $G \leq \Aut(\G)$ and a $G$-invariant partition of the edge set of $\G$ into a $2$-factor $\cC$ consisting of $8$-cycles and a $1$-factor such that $\G$ is of bialternating cycle quotient type with respect to $\cC$. When appropriate, we keep using the notation and terminology from the previous two sections. In particular, $m = |\cC|$ (and so $\G$ is of order $8m$).

In Subsection~\ref{subsec:n=8} we saw that there are four cases depending on the pair $(a,b)$, namely $(a,b) \in \{(1,4), (5,4), (4,1), (1,5)\}$. The two possibilities with $b = 4$ were analyzed together with the general case in the previous section. In particular, Proposition~\ref{pro:VT2} implies that for $n = 8$ and $b = 4$, the examples admitting a vertex-transitive subgroup of automorphisms preserving the natural $2$-factor $\cC$ are precisely the graphs of the form $\cXb(2m_0,8,1,4,\ell)$, with $m_0 \geq 2$ and $\ell \in \{2,6\}$, and the graphs of the form $\cXb(m,8,5,4,\ell)$, where either $m \geq 3$ is odd and $\ell \in \{3,7\}$, or $m \geq 4$ is even and $\ell \in \{2,6\}$.

From now on we focus on the two remaining possibilities. We therefore assume that $(a,b) \in \{(4,1), (1,5)\}$. Of course, everything that was established up to Subsection~\ref{subsec:n=8}, still holds in these two cases. In particular, we can assume that the vertices of $\G$ have been labeled in such a way that each of~\eqref{eq:non-links} and~\eqref{eq:links_init} holds. Note that, since $n = 8$, this uniquely determines all the edges of $\G$, except possibly the links between $V_{m-1}$ and $V_0$. However, Lemma~\ref{le:n=8} implies that $\ell = 6$, and so $\No{\vv{m-1,2(m-1)}} = \vv{0,6}$. It is now easy to see that the assumption that $G$ acts vertex-transitively and preserves $\cC$ enables us to completely determine the remaining three links between $V_{m-1}$ and $V_0$. We describe them explicitly in Construction~\ref{con:n=8type3} and Construction~\ref{con:n=8type4}, but first we show that $m$ is divisible by $3$. 

\begin{lemma}
With the assumptions and notation of this section, assume $(a,b) \in \{(4,1), (1,5)\}$. Then $m$ is divisible by $3$.
\end{lemma}

\begin{figure}[h]
\begin{center}
\subfigure
{
\begin{tikzpicture}[scale = .6]
\foreach \i in {0,1,2,3}{
\foreach \j in {0,1,2,3,4,5,6,7}{
\node[vtx, ,fill = white, inner sep = 1.5pt,] (A\i\j) at (2*\i, 1*\j) {};
}}
\begin{scope}[on background layer]
\foreach \j in {0,1,2,3,4,5,6}{
\foreach \i in {0,1,2,3}{
\draw[thick] let \n1 = {int(mod(\j+1, 8))} in (A\i\j) -- (A\i\n1);
}}
\foreach \i in {0,1,2,3}{
\draw[thick] (A\i 0) to[bend left = 15] (A\i 7);
}
\draw[thick] (A00) -- (A10);
\draw[thick] (A01) -- (A14);
\draw[thick] (A04) -- (A11);
\draw[thick] (A05) -- (A15);
\draw[thick] (A12) -- (A22);
\draw[thick] (A16) -- (A23);
\draw[thick] (A13) -- (A26);
\draw[thick] (A17) -- (A27);
\draw[thick] (A24) -- (A34);
\draw[thick] (A25) -- (A30);
\draw[thick] (A21) -- (A31);
\draw[thick] (A20) -- (A35);
\draw[opacity = .5, blue, line width = 2 pt] (A00) -- (A10);
\draw[opacity = .5, blue, line width = 2 pt] (A04) -- (A11);
\draw[opacity = .5, magenta, line width = 2 pt] (A01) -- (A14);
\draw[opacity = .5, magenta, line width = 2 pt] (A05) -- (A15);
\draw[opacity = .5, blue, line width = 2 pt] (A12) -- (A22);
\draw[opacity = .5, blue, line width = 2 pt] (A17) -- (A27);
\draw[opacity = .5, magenta, line width = 2 pt] (A13) -- (A26);
\draw[opacity = .5, magenta, line width = 2 pt] (A16) -- (A23);
\draw[opacity = .5, blue, line width = 2 pt] (A24) -- (A34);
\draw[opacity = .5, blue, line width = 2 pt] (A25) -- (A30);
\draw[opacity = .5, magenta, line width = 2 pt] (A20) -- (A35);
\draw[opacity = .5, magenta, line width = 2 pt] (A21) -- (A31);
\end{scope}
\end{tikzpicture}
}
\hskip 25mm
\subfigure
{
\begin{tikzpicture}[scale = .6]
\foreach \i in {0,1,2,3}{
\foreach \j in {0,1,2,3,4,5,6,7}{
\node[vtx, ,fill = white, inner sep = 1.5pt,] (A\i\j) at (2*\i, 1*\j) {};
}}
\begin{scope}[on background layer]
\foreach \j in {0,1,2,3,4,5,6}{
\foreach \i in {0,1,2,3}{
\draw[thick] let \n1 = {int(mod(\j+1, 8))} in (A\i\j) -- (A\i\n1);
}}
\foreach \i in {0,1,2,3}{
\draw[thick] (A\i 0) to[bend left = 15] (A\i 7);
}
\draw[thick] (A00) -- (A10);
\draw[thick] (A01) -- (A11);
\draw[thick] (A04) -- (A15);
\draw[thick] (A05) -- (A14);
\draw[thick] (A12) -- (A22);
\draw[thick] (A16) -- (A27);
\draw[thick] (A13) -- (A23);
\draw[thick] (A17) -- (A26);
\draw[thick] (A24) -- (A34);
\draw[thick] (A25) -- (A35);
\draw[thick] (A20) -- (A31);
\draw[thick] (A21) -- (A30);
\draw[opacity = .5, blue, line width = 2 pt] (A00) -- (A10);
\draw[opacity = .5, blue, line width = 2 pt] (A04) -- (A15);
\draw[opacity = .5, magenta, line width = 2 pt] (A01) -- (A11);
\draw[opacity = .5, magenta, line width = 2 pt] (A05) -- (A14);
\draw[opacity = .5, blue, line width = 2 pt] (A12) -- (A22);
\draw[opacity = .5, blue, line width = 2 pt] (A13) -- (A23);
\draw[opacity = .5, magenta, line width = 2 pt] (A16) -- (A27);
\draw[opacity = .5, magenta, line width = 2 pt] (A17) -- (A26);
\draw[opacity = .5, blue, line width = 2 pt] (A24) -- (A34);
\draw[opacity = .5, blue, line width = 2 pt] (A21) -- (A30);
\draw[opacity = .5, magenta, line width = 2 pt] (A25) -- (A35);
\draw[opacity = .5, magenta, line width = 2 pt] (A20) -- (A31);
\end{scope}
\end{tikzpicture}
}
\caption{The two cases when $n = 8$ and $(a,b) \in \{(4,1), (1,5)\}$.}
\label{fig:n=8exceptional}
\end{center}
\end{figure}
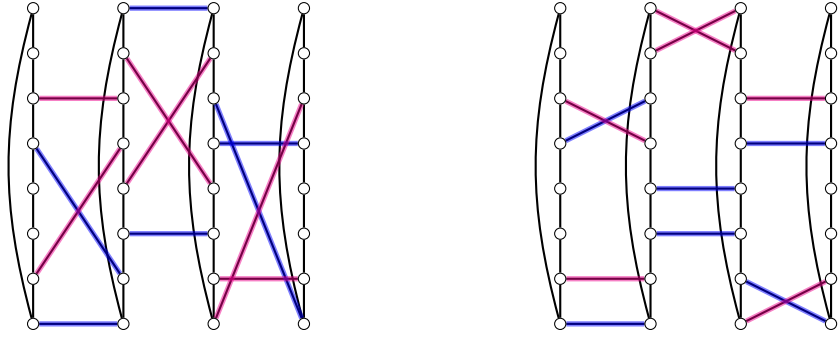

\begin{proof}
If $m = 3$, there is nothing to prove, so assume $m > 3$. By Lemma~\ref{le:G_regular}, the group $G$ is regular, and thus contains a unique automorphism $\rho$, whose action on $V_0$ is as in~\eqref{eq:rho_V_0}. Recall that the involution $\rho^2$ preserves each $V_i$, $i \in \ZZ_m$, and consider its action on each of the cycles $C_i$. For each of the two possibilities, the part of the graph $\G$, containing the vertices from the union $V_0 \cup V_1 \cup V_2 \cup V_3$, is presented on Figure~\ref{fig:n=8exceptional} (the case $(a,b) = (4,1)$ on the left and the case $(a,b) = (1,5)$ on the right part of the figure). On this figure, the two links connecting $V_i$ to $V_{i+1}$, $i \in \{0,1,2\}$, that are interchanged by the involution $\rho^2$, are highlighted with the same color. It is now clear that the involution $\rho^2$ is a (step $4$) rotation on $C_0$, is not a rotation on any of $C_1$ and $C_2$, and is again a (step $4$) rotation on $C_3$. An inductive argument thus shows that $\rho^2$ is a rotation on each $C_i$ with $3 \mid i$, and is not a rotation on any other $C_i$. Since $\G$ is $G$-vertex-transitive, the same pattern (of a rotational action following two consecutive non-rotational actions) thus must also continue when we return from $C_{m-1}$ back to $C_0$. Consequently, $m$ is divisible by $3$, as claimed.
\end{proof}

In the next two subsections we show that for each $m$ with $3 \mid m$, each of the graphs corresponding to $(a,b) \in \{(4,1), (1,5)\}$ (and where $\No{\vv{m-1,2(m-1)}} = \vv{0,6}$), does indeed admit a regular group preserving the natural $2$-factor $\cC$. 

\subsection{The $(a,b) = (4,1)$ case}

Consider first the case that $(a,b) = (4,1)$. Based on what was established so far in this section we only need to determine the remaining three links between the vertices of $V_{m-1}$ and $V_0$. Since $\No{\vv{0,0}}$ and $\No{\vv{0,1}}$ are antipodal on $C_1$, so are $\vv{0,6}$ ($ = \No{\vv{m-1,2(m-1)}}$) and $\No{\vv{m-1,2(m-1)+1}}$ on $C_0$, and so $\No{\vv{m-1,2(m-1)+1}} = \vv{0,2}$. Likewise, $\vv{m-1,2(m-1)}$ and $\No{\vv{0,7}}$ must be antipodal on $C_{m-1}$, and so $\No{\vv{m-1,2(m-1)+4}} = \vv{0,7}$. It is now clear that our graphs in the case of $(a,b) = (4,1)$ belong to the following infinite family.

\begin{construction}
\label{con:n=8type3}
Let $m \geq 3$ be an integer divisible by $3$. Then $\cXb^1(m)$ is the graph with vertex set $\{\vv{i,j} \colon i \in \ZZ_{m}, j \in \ZZ_8\}$ in which $\vv{i,j} \sim \vv{i,j+1}$ for all $i \in \ZZ_m$, $j \in \ZZ_8$, 
$$
	\vv{i,2i+\delta} \sim \vv{i+1,2i+4\delta}\ \text{and}\ \vv{i,2i+4+\delta} \sim \vv{i+1,2i+1+4\delta}
$$
for all $i \in \{0,1,\ldots , m-2\}$ and $\delta \in \{0,1\}$, and 
$$
	\vv{m-1,2(m-1)+\delta} \sim \vv{0,6+4\delta}\ \text{and}\ \vv{m-1,2(m-1)+4+\delta} \sim \vv{0,7+4\delta}
$$
for $\delta \in \{0,1\}$.
\end{construction}

We let the meaning of the sets $V_i$, cycles $C_i$, the natural $2$-factor $\cC$, and the links and non-links of a graph of the form~$\cXb^1(m)$ be consistent with what was agreed for the graphs $\cXb(m,n,a,b,\ell)$ in the paragraph following Construction~\ref{con:Xb}.

\begin{proposition}
\label{pro:n=8_T3}
Let $m \geq 3$ be an integer divisible by $3$ and let $\G = \cXb^1(m)$. Then $\G$ admits a regular subgroup $G \leq \Aut(\G)$, preserving the natural $2$-factor $\cC$ consisting of $8$-cycles. However, the full automorphism group $\Aut(\G)$ has vertex stabilizers of order $2$, and so $\cC$ is not $\Aut(\G)$-invariant.
\end{proposition}

\begin{figure}[h]
\begin{center}
\subfigure
{
\begin{tikzpicture}[scale = .6]
\foreach \i in {0,1,2,3}{
\foreach \j in {0,1,2,3,4,5}{
\node[vtx, ,fill = white, inner sep = 1.5pt,] (A\i\j) at (360*\j/6 : \i+2) {};
}}
\begin{scope}[on background layer]
\node[right = -1pt of A30] {\footnotesize $\vv{0,0}$};
\node[right = -1pt of A31] {\footnotesize $\vv{1,0}$};
\node[left = -1pt of A32] {\footnotesize $\vv{1,7}$};
\node[left = -1pt of A33] {\footnotesize $\vv{2,7}$};
\node[left = -1pt of A34] {\footnotesize $\vv{2,0}$};
\node[right = -1pt of A35] {\footnotesize $\vv{0,7}$};
\node[left = -1pt of A00] {\footnotesize $\vv{0,5}$};
\node[right = -1pt of A03] {\footnotesize $\vv{2,3}$};
\node[below left = -2pt of A01] {\footnotesize $\vv{1,5}$};
\node[above right = -3pt of A04] {\footnotesize $\vv{2,4}$};
\node[above  left = -3pt of A21] {\footnotesize $\vv{1,4}$};
\node[above  right = -3pt of A22] {\footnotesize $\vv{1,3}$};
\foreach \i in {0,1,2,3}{
\foreach \j in {0,1,2,3,4,5}{
\draw[thick] let \n1 = {int(mod(\j+1, 6))} in (A\i\j) -- (A\i\n1);
}}
\draw[thick] (A00) -- (A10);
\draw[thick] (A03) -- (A13);
\draw[thick] (A20) -- (A30);
\draw[thick] (A23) -- (A33);
\draw[thick] (A01) to[bend left = 30] (A21);
\draw[thick] (A04) to[bend left = 30] (A24);
\draw[thick] (A11) to[bend right = 30] (A31);
\draw[thick] (A14) to[bend right = 30] (A34);
\draw[thick] (A02) to[bend left = 30] (A32);
\draw[thick] (A05) to[bend left = 30] (A35);
\draw[thick] (A12) -- (A22);
\draw[thick] (A15) -- (A25);
\end{scope}
\end{tikzpicture}
}
\hskip 5mm
\subfigure
{
\begin{tikzpicture}[scale = .6]
\foreach \i in {0,1,2,3}{
\foreach \j in {0,1,2,3,4,5,6,7,8,9,10,11}{
\node[vtx, fill = white, inner sep = 1.5pt,] (A\i\j) at (360*\j/12 : \i+2) {};
}}
\begin{scope}[on background layer]
\node[right = -1pt of A30] {\footnotesize $\vv{0,0}$};
\node[right = -1pt of A31] {\footnotesize $\vv{1,0}$};
\node[right = -1pt of A32] {\footnotesize $\vv{1,7}$};
\node[left = -1pt of A34] {\footnotesize $\vv{2,0}$};
\node[left = -1pt of A35] {\footnotesize $\vv{3,5}$};
\node[left = -1pt of A37] {\footnotesize $\vv{4,6}$};
\node[left = -1pt of A38] {\footnotesize $\vv{4,5}$};
\node[right = -1pt of A310] {\footnotesize $\vv{5,6}$};
\node[right = -1pt of A311] {\footnotesize $\vv{0,7}$};
\node[left = -1pt of A00] {\footnotesize $\vv{0,5}$};
\node[left = -1pt of A01] {\footnotesize $\vv{1,5}$};
\node[left = -1pt of A011] {\footnotesize $\vv{0,6}$};
\node[below = -1pt of A03] {\footnotesize $\vv{2,3}$};
\node[right = -1pt of A06] {\footnotesize $\vv{3,3}$};
\node[above = -1pt of A09] {\footnotesize $\vv{5,1}$};
\foreach \i in {0,1,2,3}{
\foreach \j in {0,1,2,3,4,5,6,7,8,9,10,11}{
\draw[thick] let \n1 = {int(mod(\j+1, 12))} in (A\i\j) -- (A\i\n1);
}}
\draw[thick, bend left = 30] (A01) to (A21);
\draw[thick, bend right = 30] (A11) to (A31);
\draw[thick, bend left = 30] (A04) to (A24);
\draw[thick, bend right = 30] (A14) to (A34);
\draw[thick, bend left = 30] (A07) to (A27);
\draw[thick, bend right = 30] (A17) to (A37);
\draw[thick, bend left = 30] (A010) to (A210);
\draw[thick, bend right = 30] (A110) to (A310);
\draw[thick, bend left = 30] (A02) to (A32);
\draw[thick] (A12) -- (A22);
\draw[thick, bend left = 30] (A05) to (A35);
\draw[thick] (A15) -- (A25);
\draw[thick, bend left = 30] (A08) to (A38);
\draw[thick] (A18) -- (A28);
\draw[thick, bend left = 30] (A011) to (A311);
\draw[thick] (A111) -- (A211);
\foreach \j in {0,3,6,9}{
\draw[thick] (A0\j) -- (A1\j);
\draw[thick] (A2\j) -- (A3\j);
}
\end{scope}
\end{tikzpicture}
}
\caption{A different presentation of the graphs $\cXb^1(3)$ and $\cXb^1(6)$.}
\label{fig:cXb1}
\end{center}
\end{figure}
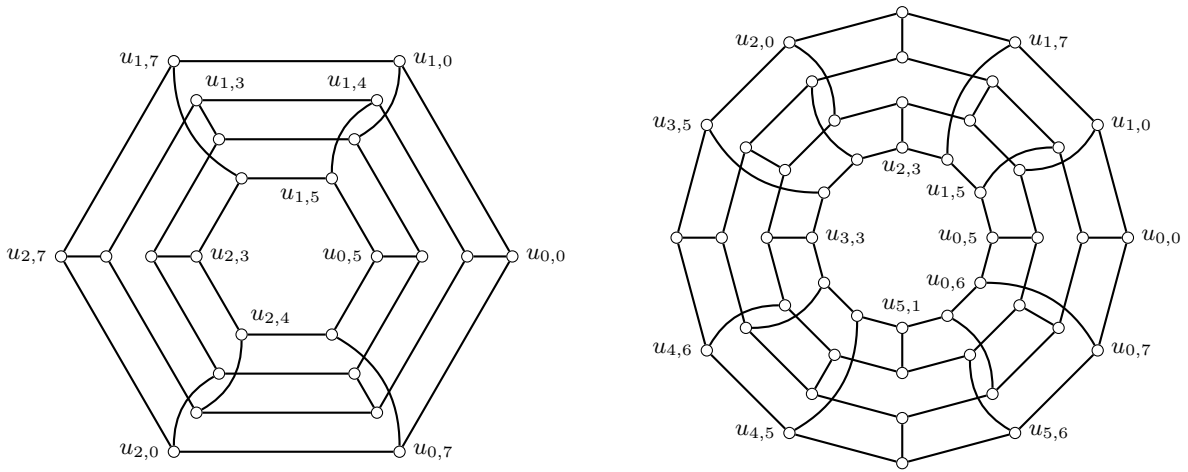

\begin{proof}
We point out the main steps of the proof but leave out some of the rather tedious verifications that the stated permutations are indeed automorphisms of $\G$. 

It is easy to verify that the permutation $\gamma$, defined by~\eqref{eq:gamma} and~\eqref{eq:gamma_m1} (where we take $\ell = 6$), is an automorphism of $\G$ (what really needs to be checked are only the links between $V_{m-2}$ and $V_{m-1}$ and between $V_{m-1}$ and $V_0$). This shows that the subgraphs $\G_i$, $i \in \ZZ_m$, of $\G$, induced on the union $V_i \cup V_{i+1} \cup V_{i+2} \cup V_{i+3}$ (where we compute the indices modulo $m$), are all isomorphic. Figure~\ref{fig:n=8exceptional} reveals that there exists an automorphism of $\G_0$, acting as a $4$-step rotation on each of $C_0$ and $C_{3}$, and as a reflection on each of $C_{1}$ and $C_{2}$. It thus follows that there exist $\varphi_0, \varphi_1 \in \Aut(\G)$, each of which preserves each $V_i$, $i \in \ZZ_m$, and such that $\varphi_0$ acts as a $4$-step rotation on each $C_i$ with $3 \mid i$, while $\varphi_1$ acts as a $4$-step rotation on each $C_i$ with $3 \mid (i-1)$. Since $\varphi_1$ interchanges $\vv{0,0}$ and $\vv{0,1}$, while $\varphi_0$ maps these two vertices to $\vv{0,4}$ and $\vv{0,5}$, respectively, it clearly follows that the group $\la \gamma, \varphi_0, \varphi_1 \ra$ has two orbits on the vertex set of $\G$.

Write $m = 3m_0$ and let $\alpha$ be the permutation of the vertex set of $\G$, where for each $j \in \ZZ_8$ we set the following: $\alpha$ interchanges $\vv{0,j}$ with $\vv{1,j}$, while for each $i$ with $1 \leq i \leq m_0$ and $\varepsilon \in \{0,1,2\}$ such that $1 \leq 3i-\varepsilon \leq m-2$, we set
\[
	\alpha(\vv{m-3i+\varepsilon, j}) = \left\{\begin{array}{ccc} 
		\vv{3i+1-\varepsilon, 2m_0+4i+j} & ; & \varepsilon \in \{0,1\} \\
		\vv{3i-1, 1-2m_0-j} & ; & \varepsilon = 2.
		\end{array}\right.
\]
A somewhat tedious but fairly straightforward computation shows that $\alpha$ is an automorphism of $\G$ (clearly, we only need to verify that the links are mapped to links of $\G$). For instance, the vertices $\vv{1,2}$, $\vv{1,3}$, $\vv{1,6}$ and $\vv{1,7}$ are mapped by $\alpha$ to $\vv{0,2}$, $\vv{0,3}$, $\vv{0,6}$ and $\vv{0,7}$, respectively. Since $2 = m - 3m_0 + 2$ and $2(m-1) = -2m_0 - 2$ in $\ZZ_8$, their corresponding outside neighbors $\vv{2,2}$, $\vv{2,6}$, $\vv{2,3}$ and $\vv{2,7}$, respectively, are mapped to $\vv{m-1,2(m-1)+1}$, $\vv{m-1,2(m-1)+5}$, $\vv{m-1,2(m-1)}$ and $\vv{m-1,2(m-1)+4}$, respectively. This confirms that the links between $V_1$ and $V_2$ are indeed mapped to links of $\G$. We leave the rest of the verification to the reader. Since the group $\la \gamma, \varphi_0, \varphi_1, \alpha \ra$ clearly preserves the $2$-factor $\cC$ and acts transitively on the vertex set of $\G$, the first part of the proof is complete by Lemma~\ref{le:G_regular}. 

To prove the second part, observe first that $\G$ possesses $7$-cycles. It is clear that if $m > 3$, then each $7$-cycle of $\G$ is contained in a subgraph of $\G$, induced on $V_i \cup V_{i+1}$ for some $i \in \ZZ_m$. One can check that the same actually holds for the graph $\cXb^1(3)$. It is now not difficult to verify that the edge $\vv{0,0}\vv{0,1}$ lies on $6$ different $7$-cycles of $\G$, while each of the edges $\vv{0,0}\vv{1,0}$ and $\vv{0,0}\vv{0,7}$ lies on $4$ different $7$-cycles of $\G$. This implies that the set $\{\vv{i,j}\vv{i,j+1} \colon i \in \ZZ_m,\ j \in \{0,2,4,6\}\}$ is an $\Aut(\G)$-orbit. Consequently, an automorphism of $\G$, fixing the vertex $\vv{0,0}$ and each of its three neighbors $\vv{1,0}$, $\vv{0,1}$ and $\vv{0,7}$, necessarily also fixes each neighbor of $\vv{1,0}$ and each neighbor of $\vv{0,7}$. Since there is precisely one $7$-cycle containing the $3$-path $(\vv{0,1},\vv{0,0},\vv{1,0},\vv{1,7})$, we see that such an automorphism also fixes each neighbor of $\vv{0,1}$. Since $\G$ is connected and vertex-transitive, this finally shows that an automorphism of $\G$, fixing a vertex and each of its three neighbors, is necessarily the identity, thus proving that $\Aut(\G)$ has vertex stabilizers of order at most $2$. 

It turns out however, that a nontrivial automorphism of $\G$, fixing $\vv{0,0}$ (and interchanging $\vv{1,0}$ with $\vv{0,7}$) indeed exists. Instead of spelling out the action of this automorphism and proving that it preserves adjacency, we make the following observation. We claim that if we delete the $1$-factor of $\G$ corresponding to the above mentioned $\Aut(\G)$-orbit of the edge $\vv{0,0}\vv{0,1}$, we are left with a $2$-factor of $\G$, consisting of four $(2m)$-cycles. To see this let us follow the cycle starting with the vertices $\vv{0,0},\vv{1,0},\vv{1,7},\vv{2,7} (= \vv{2,2\cdot 2 + 3}),\vv{2,0} (= \vv{2,2\cdot 2 + 4})$. We see that for each $i \in \ZZ_m$ with $3 \mid i$, we enter the set $V_{i-1}$ in the vertex $\vv{i-1,2(i-1)+3}$ and then leave it from the vertex $\vv{i-1,2(i-1)+4}$. Consequently, the cycle is indeed of length $2m$ and ``ends'' with the three vertices $\vv{m-1,2(m-1)+3},\vv{m-1,2(m-1)+4}$ and $\vv{0,7}$ (recall that $3 \mid m$). With this in mind, the examples $\cXb^1(3)$ and $\cXb^1(6)$ are presented on Figure~\ref{fig:cXb1} in a way from which the existence of the above mentioned automorphism (in a general graph $\cXb^1(m)$ with $3 \mid m$) is easy to justify - we reflect with respect to the line through $\vv{0,0}$ and $\vv{0,5}$ and then swap the two ``inner'' $(2m)$-cycles.
\end{proof}

\subsection{The $(a,b) = (1,5)$ case}

We finally consider the case that $(a,b) = (1,5)$. In a similar way as was done in the case of $(a,b) = (4,1)$ we can show that this time we need to consider the graphs from the following infinite family.

\begin{construction}
\label{con:n=8type4}
Let $m \geq 3$ be an integer divisible by $3$. Then $\cXb^2(m)$ is the graph with vertex set $\{\vv{i,j} \colon i \in \ZZ_{m}, j \in \ZZ_8\}$ in which $\vv{i,j} \sim \vv{i,j+1}$ for all $i \in \ZZ_m$, $j \in \ZZ_8$, 
$$
	\vv{i,2i+\delta} \sim \vv{i+1,2i+\delta}\ \text{and}\ \vv{i,2i+4+\delta} \sim \vv{i+1,2i+5-\delta}
$$
for all $i \in \{0,1,\ldots , m-2\}$ and $\delta \in \{0,1\}$, and 
$$
	\vv{m-1,2(m-1)+\delta} \sim \vv{0,6+\delta}\ \text{and}\ \vv{m-1,2(m-1)+4+\delta} \sim \vv{0,3-\delta}
$$
for $\delta \in \{0,1\}$.
\end{construction}

Again, we let the meaning of the sets $V_i$, cycles $C_i$, the natural $2$-factor $\cC$, and the links and non-links of a graph of the form~$\cXb^2(m)$ be consistent with what was agreed for the graphs $\cXb(m,n,a,b,\ell)$ in the paragraph following Construction~\ref{con:Xb}.

\begin{proposition}
\label{pro:n=8_T4}
Let $m \geq 3$ be an integer divisible by $3$ and let $\G = \cXb^2(m)$. Then the automorphism group of $\G$ acts regularly on its vertex set and preserves the natural $2$-factor $\cC$ consisting of $8$-cycles. 
\end{proposition}

\begin{proof}
The proof is very similar to that of Proposition~\ref{pro:n=8_T3}. We again first verify that the permutation $\gamma$, defined by~\eqref{eq:gamma} and~\eqref{eq:gamma_m1} (where we take $\ell = 6$), is an automorphism of $\G$. We can also easily justify the existence of $\varphi_0, \varphi_1 \in \Aut(\G)$, preserving each $V_i$, $i \in \ZZ_m$, and such that $\varphi_0$ acts as a $4$-step rotation on $V_i$ for each $i \in \ZZ_m$ with $3 \mid i$ and as a reflection on all other $V_i$, and that $\varphi_1$ acts as a $4$-step rotation on $V_i$ for each $i \in \ZZ_m$ with $3 \mid (i-1)$ and as a reflection on all other $V_i$. To conclude the proof we thus only need to confirm the existence of an automorphism $\alpha$, preserving the $2$-factor $\cC$ and interchanging $\vv{0,0}$ with $\vv{1,0}$, and to verify that the identity is the only automorphism of $\G$ fixing $\vv{0,0}$.

For the first part we write $m = 3m_0$ and let $\alpha$ be the permutation of the vertex set of $\G$, where for each $j \in \ZZ_8$ we set the following: $\alpha$ interchanges $\vv{0,j}$ with $\vv{1,j}$, while for each $i$ with $1 \leq i \leq m_0$ and $\varepsilon \in \{0,1,2\}$ such that $1 \leq 3i-\varepsilon \leq m-2$, we set
\[
	\alpha(\vv{m-3i+\varepsilon, j}) = \left\{\begin{array}{ccc} 
		\vv{3i+1-\varepsilon, 2m_0+4i+j} & ; & \varepsilon \in \{0,1\} \\
		\vv{3i-1, 5-2m_0-j} & ; & \varepsilon = 2.
		\end{array}\right.
\]
A straightforward calculation shows that $\alpha$ preserves the set of the links of $\G$, and is thus an automorphism of $\G$. Lemma~\ref{le:G_regular} then implies that the group $\la \gamma, \varphi_0, \varphi_1, \alpha \ra$ acts regularly on $\G$. To see that in fact $\Aut(\G) = \la \gamma, \varphi_0, \varphi_1, \alpha \ra$, we can do the following. We first observe that each of the edges $\vv{0,0}\vv{0,1}$ and $\vv{0,0}\vv{1,0}$ lies on a $4$-cycle of $\G$, while $\vv{0,0}\vv{0,7}$ does not. Therefore, $\Aut(\G)$ is regular on the vertex set of $\G$ if and only if there is no automorphism of $\G$, fixing $\vv{0,0}$ but swapping $\vv{0,1}$ with $\vv{1,0}$. If $m = 3$, then one readily verifies that the edge $\vv{0,0}\vv{1,0}$ lies on a $6$-cycle of $\G$, while $\vv{0,0}\vv{0,1}$ does not. If however, $m > 3$, then clearly no $8$-cycle of $\G$ can contain at least one vertex from each of the sets $V_i$, and it is then easy to conclude that the only $8$-cycles of $\G$ are the members of the $2$-factor $\cC$. This clearly shows that no automorphism of $\G$ can map the edge $\vv{0,0}\vv{0,1}$ to the edge $\vv{0,0}\vv{1,0}$, completing the proof.
\end{proof}

We can now state the main result of this paper. Before doing so we point out that since the parameter $a$ is of the form $a = 4a_0 + 1$, the walk $(\vv{2,3},\vv{2,2},\vv{1,2},\vv{1,1},\vv{1,0},\vv{0,0},\vv{0,1},\vv{1,a},\vv{1,a+1},\vv{1,a+2})$ in a graph of the form $\cXb(m,n,a,b,\ell)$ is a path (of length $9$) if and only if $a \neq 1$. Therefore, since~\eqref{eq:b_a_cond} and~\eqref{eq:links} imply that $\No{\vv{1,a+2}} = \vv{2,2+ba_0+a} = \vv{2,3}$, this shows that if $a \neq 1$ the above path in fact closes to a $10$-cycle of $\G$. This proves the claim regarding the girth in the following theorem.

\begin{theorem}
\label{the:main}
Let $\G$ be a cubic vertex-transitive graph admitting a partition of its edge set into a $2$-factor $\cC$ and a $1$-factor such that $\G$ is of bialternating cycle quotient type with respect to $\cC$. Let $m = |\cC|$ and let $n$ be the length of the cycles from $\cC$. Then $\G$ admits a vertex-transitive subgroup $G$ of $\Aut(\G)$ preserving the $2$-factor $\cC$ if and only if $4 \mid n$ and one of the following holds:
\begin{enumerate}\itemsep = 0pt
\item[(i)] $n = 4$ and $\G$ is isomorphic to the M\"obius ladder or the prism graph of order $4m$;
\item[(ii)] $n = 8$ and one of the following holds:
	\begin{itemize}\itemsep = 0pt
	\item $3 \mid m$ and $\G \cong \cXb^1(m)$ or $\G \cong \cXb^2(m)$;
	\item $m$ is odd and $\G \cong \cXb(m,8,5,4,\ell)$ for $\ell \in \{3,7\}$;
	\item $m$ is even and $\G \cong \cXb(m,8,a,4,\ell)$ for some $a \in \{1,5\}$ and $\ell \in \{2,6\}$.
	\end{itemize}
\item[(iii)] $n \geq 12$ and $\G \cong \cXb(m,n,a,b,\ell)$, where one of the following holds:
	\begin{itemize}\itemsep = 0pt
	\item $m$ and $\ell$ are both odd, $b = n-4$ and $a = 2\ell-1$ holds in $\ZZ_n$;
	\item $m$ is odd, $\ell$ is even, $a = 1$, $b = 4$, $n = 4n_0$ for an odd integer $n_0 \geq 3$, and $\ell = 2n_0 + 2m - 2$ holds in $\ZZ_n$.
	\item $m$ and $\ell$ are both even, and $(\ell + 2)(b_0 - 1) = 0$ and $2(\ell + 2) = 2m(b_0+1)$ both hold in $\ZZ_n$, where $b = 4b_0$.
	\end{itemize}
\end{enumerate}
Moreover, $\G$ is of girth at most $10$, the group $G$ is regular, and $\G$ is a Cayley graph of $G$ with respect to a connection set consisting of three involutions.
\end{theorem}

\section{Additional automorphisms}
\label{sec:addAut}

Theorem~\ref{the:main} provides a complete solution of~\cite[Problem 1.1]{AlsSpa24} in the case of the graphs of bialternating cycle quotient type. However, for some of the graphs $\G$ from Theorem~\ref{the:main}, the natural $2$-factor is not $\Aut(\G)$-invariant. Recall that by Lemma~\ref{le:G_regular}, this happens if and only if the group $\Aut(\G)$ has nontrivial vertex-stabilizers. It turns out that, at least for the examples from item (iii) of the theorem, determining when this can happen is not an easy task. In this paper we settle this question for the graphs from items (i) and (ii) of the theorem, indicate how this can be done for the second of the three possibilities from item (iii) of the theorem, and leave the analysis of the remaining graphs $\cXb(m,n,a,b,\ell)$ with $n \geq 12$ for future investigations. 
\medskip

We already mentioned that the graphs from item (i), that is the M\"obius ladders and the prism graphs, do not have regular automorphism groups. In the case of the graphs from item (i) the natural $2$-factor $\cC$ (consisting of $4$-cycles) is therefore not $\Aut(\G)$-invariant. 
\medskip

The graphs from item (ii) of the theorem (those with $n = 8$) require a bit more work, but using some known results, the task is not difficult. Proposition~\ref{pro:n=8_T3} and Proposition~\ref{pro:n=8_T4} show that the natural $2$-factor $\cC$ is not $\Aut(\G)$-invariant for the graphs $\G = \cXb^1(m)$ with $3 \mid m$, but is $\Aut(\G)$-invariant for the graphs $\G = \cXb^2(m)$ with $3 \mid m$. 

For the remaining examples with $n = 8$ we indicate the main idea of the analysis, but leave the details to the reader. For the examples of the form $\G = \cXb(m,8,5,4,\ell)$, the key idea is to verify that these graphs are in fact isomorphic to certain honeycomb toroidal graphs (for the definition of the $\mathrm{HTG}(m,n,\ell)$ graphs we refer the reader to~\cite{Spa22}). To see this, we consider the spanning subgraph of $\G$, given by all the links and by all the non-links of $\G$ of the form $\vv{i,j}\vv{i,j+1}$ with $j$ odd. It is not difficult to verify that this subgraph consists of two $(4m)$-cycles when either $m$ is odd and $\ell = 3$ or $m$ is even and $\ell = 2$, while it consists of four $(2m)$-cycles when either $m$ is odd and $\ell = 7$ or $m$ is even and $\ell = 6$. It is then easy to confirm that in the case when $m$ is odd, we obtain $\cXb(m,8,5,4,3) \cong \mathrm{HTG}(2,4m,2m)$ and $\cXb(m,8,5,4,7) \cong \mathrm{HTG}(4,2m,0)$, while in the case when $m$ is even, we obtain $\cXb(m,8,5,4,2) \cong \mathrm{HTG}(2,4m,2m)$ and $\cXb(m,8,5,4,6) \cong \mathrm{HTG}(4,2m,0)$. Then~\cite[Theorem~1.2]{Spa22} implies that $\cXb(3,8,5,4,3) \cong \mathrm{HTG}(2,12,6)$ is in fact $2$-arc-regular, while all the other graphs have vertex-stabilizers of order $2$. Consequently, the natural $2$-factor $\cC$ is not $\Aut(\G)$-invariant in any of these graphs.

To complete the $n = 8$ case, we thus finally need to consider the examples $\cXb(m,8,1,4,\ell)$ with $m \geq 4$ even and $\ell \in \{2,6\}$. That the graph $\cXb(4,8,1,4,6)$ has vertex-stabilizers of order $4$ (and so $\cC$ is not $\Aut(\G)$-invariant), can be verified using a computer. That all other graphs of the form $\cXb(m,8,1,4,\ell)$ with $m \geq 4$ even and $\ell \in \{2,6\}$ have regular automorphism groups, follows from the following observation. Of the three edges incident to $\vv{1,2}$, only $\vv{1,1}\vv{1,2}$ does not lie on a $4$-cycle, while of the remaining two edges the edge $\vv{1,2}\vv{2,2}$ lies on a unique $8$-cycle, whereas the edge $\vv{1,2}\vv{1,3}$ lies on two. This shows that the three mentioned edges are in three different $\Aut(\G)$-orbits, and so the fact that $\Aut(\G)$ is vertex-transitive, implies that it is actually regular.
\medskip

As we already mentioned, determining which of the graphs from item (iii) of Theorem~\ref{the:main} have regular automorphism groups seems to be rather difficult. The main reason is that, except in the case that $a = 1$, the graphs can have girth $10$, which makes the usual approach via the analysis of short cycles rather unappealing. 

At this point we therefore content ourselves by only completing the analysis of the examples with $m$ odd and $\ell$ even. Recall that Theorem~\ref{the:main} implies that in this case $a = 1$, and so $\G$ has $4$-cycles. Moreover, $b = 4$, $n = 4n_0$ for an odd $n_0 \geq 3$, and $\ell = 2n_0 + 2m-2$ holds in $\ZZ_n$. Let thus $\G = \cXb(m,n,1,4,\ell)$ for an odd $m \geq 3$ and such $n$ and $\ell$. Note that the non-links of the form $\vv{i,j}\vv{i,j+1}$, where $j$ is odd, do not lie on a $4$-cycle, while all the other edges of $\G$ do. Therefore, $\mathcal{G} = \{\{\vv{i,j}\vv{i,j+1}\} \colon i \in \ZZ_m, j \in \ZZ_n, 2 \nmid j\}$ is an $\Aut(\G)$-orbit. We think of coloring these edges green, coloring all the links red and coloring the remaining non-links blue. Note that the three color classes coincide with the three $G$-orbits on the edge set of $\G$, where $G$ is the regular group from Theorem~\ref{the:main}.

The fact that each blue and each red edge lies on a unique $4$-cycle, while the green edges do not lie on $4$-cycles, imposes severe restrictions on the possibilities for automorphisms outside the regular group $G$. One could proceed with the analysis as follows. It is easy to see that each $4$-cycle of $\G$ is ``adjacent'' to four different $4$-cycles, and so with the terminology of~\cite{EibJajSpa19}, the graph $\G$ is a generalized truncation of an arc-transitive $4$-valent graph $\Delta$ of girth $4$. One could now use the results of~\cite{EibJajSpa19}, which show that the full automorphism group of $\G$ ``projects'' to $\Delta$, then use the results of~\cite{PotWil07} to see that $\Delta$ must be the skeleton of an edge-transitive map of type $\{4,4\}$ on the torus, and finally use the classification of these maps~\cite{SirTucWat01} to search for those that admit additional automorphisms. But one can also use the following simple direct argument, for which we only give the main idea.

Since $\mathcal{G}$ is an $\Aut(\G)$-orbit, an automorphism $\eta$ of $\G$, fixing the vertex $\vv{1,2}$ and each of its three neighbors, fixes each vertex at distance $2$ from $\vv{1,2}$, except possibly interchanging $\vv{0,1}$ with $\vv{1,0}$. 
However, if $\eta$ does interchange $\vv{0,1}$ with $\vv{1,0}$, it also interchanges $\vv{0,2}$ with $\vv{1,-1}$. But then the fact that there exists a path of length $3$ from $\vv{0,2}$ to the fixed vertex $\vv{1,4} = \No{\vv{0,4}}$, implies that there is also a path of length $3$ from $\vv{1,-1}$ to $\vv{1,4}$. Since $n \geq 12$, this is impossible, proving that $\eta$ fixes each vertex at distance $2$ from $\vv{1,2}$, and so connectedness and vertex-transitivity of $\G$ implies that in fact $\eta = 1$. This shows that either $G = \Aut(\G)$ or $G$ is an index $2$ subgroup of $\Aut(\G)$. 

To complete the analysis assume the latter case occurs and note that in this case each automorphism outside $G$ interchanges the set of the red edges with the set of the blue ones. Let $\eta \in \Aut(\G)$ be the unique automorphism fixing $\vv{1,1}$ but interchanging $\vv{0,1}$ with $\vv{1,0}$. It then maps the $n$-cycle $C_1$ to the cycle containing the edge $\vv{0,1}\vv{1,1}$ on which the edges alternate between being red and green. The first few vertices of this cycle, starting at $\vv{0,1}$ and going in the direction of $\vv{1,1}$, are 
$$
\vv{0,1},\vv{1,1},\vv{1,2},\vv{2,2},\vv{2,1},\vv{3,1},\vv{3,2},\vv{4,2},\vv{4,1},\ldots
$$
It thus clearly follows that $n$ must be an even multiple of $2m$ (recall that $m$ is odd and $n = 4n_0$). Moreover, $\eta$ maps the vertices $\vv{1,2(m-1)}$, $\vv{1,2m-1}$ and $\vv{1,2m}$ to $\vv{m-1,1}$, $\No{\vv{m-1,1}} = \vv{0,2n_0+1}$ and $\vv{0,2n_0+2}$, respectively. Continuing in this way we than also find that $\eta(\vv{1,4m-2}) = \vv{m-1,2+2n_0}$, $\eta(\vv{1,4m-1}) = \vv{0,2}$ and $\eta(\vv{1,4m}) = \vv{0,1} = \eta({\vv{1,0}})$, finally yielding $n = 4m$, and consequently $\ell = -2$. 

\begin{figure}[h]
\begin{center}
\subfigure
{
\begin{tikzpicture}[scale = .4]
\xdef\ssquare{1.4}
\foreach \i in {0,1,2,3,4}{
\pgfmathtruncatemacro{\xx}{\i*5};
\foreach \j in {0,1,2,3,4}{
\pgfmathtruncatemacro{\yy}{\j*5};
\node[vtx, ,fill = black, inner sep = 1.5pt,] (A\i\j0) at (\xx, \yy + \ssquare) {};
\node[vtx, ,fill = black, inner sep = 1.5pt,] (A\i\j1) at (\xx - \ssquare, \yy) {};
\node[vtx, ,fill = black, inner sep = 1.5pt,] (A\i\j2) at (\xx, \yy - \ssquare) {};
\node[vtx, ,fill = black, inner sep = 1.5pt,] (A\i\j3) at (\xx + \ssquare, \yy) {};
}}
\node[above left = -5pt of A001] {\footnotesize $u_{0,0}$};
\node[below right = -4pt of A003] {\footnotesize $u_{1,1}$};
\node[right = -3pt of A000] {\footnotesize $u_{0,1}$};
\node[above left = -5pt of A101] {\footnotesize $u_{1,2}$};
\node[below right = -4pt of A103] {\footnotesize $u_{2,3}$};
\node[right = -3pt of A012] {\footnotesize $u_{0,2}$};
\node[right = -3pt of A112] {\footnotesize $u_{1,4}$};
\node[below right = -4pt of A013] {\footnotesize $u_{0,3}$};
\node[above left = -5pt of A111] {\footnotesize $u_{0,4}$};
\node[right = -3pt of A110] {\footnotesize $u_{0,5}$};
\node[right = -3pt of A100] {\footnotesize $u_{1,3}$};
\node[below right = -4pt of A113] {\footnotesize $u_{1,5}$};
\node[right = -3pt of A122] {\footnotesize $u_{0,6}$};
\node[below right = -4pt of A123] {\footnotesize $u_{0,7}$};
\node[above left = -5pt of A221] {\footnotesize $u_{0,8}$};
\node[right = -3pt of A220] {\footnotesize $u_{0,9}$};
\node[right = -3pt of A232] {\footnotesize $u_{0,10}$};
\node[below right = -4pt of A233] {\footnotesize $u_{0,11}$};
\node[above left = -5pt of A331] {\footnotesize $u_{0,12}$};
\node[right = -3pt of A330] {\footnotesize $u_{0,13}$};
\node[right = -3pt of A342] {\footnotesize $u_{0,14}$};
\node[below right = -4pt of A343] {\footnotesize $u_{0,15}$};
\node[above left = -5pt of A441] {\footnotesize $u_{0,16}$};
\node[right = -3pt of A440] {\footnotesize $u_{0,17}$};
\node[right = -3pt of A402] {\footnotesize $u_{0,18}$};
\node[below right = -4pt of A403] {\footnotesize $u_{0,19}$};
\node[above left = -5pt of A201] {\footnotesize $u_{2,4}$};
\node[above left = -5pt of A301] {\footnotesize $u_{3,6}$};
\node[above left = -5pt of A401] {\footnotesize $u_{4,8}$};
\node[right = -3pt of A200] {\footnotesize $u_{2,5}$};
\node[below right = -4pt of A203] {\footnotesize $u_{3,5}$};
\node[right = -3pt of A300] {\footnotesize $u_{3,7}$};
\node[below right = -4pt of A303] {\footnotesize $u_{4,7}$};
\node[right = -3pt of A400] {\footnotesize $u_{4,9}$};
\node[right = -3pt of A002] {\footnotesize $u_{1,0}$};
\begin{scope}[on background layer]
\foreach \i in {0,1,2,3,4}{
\foreach \j in {0,1,2,3,4}{
\draw[thick] (A\i\j0) -- (A\i\j1);
\draw[thick] (A\i\j1) -- (A\i\j2);
\draw[thick] (A\i\j2) -- (A\i\j3);
\draw[thick] (A\i\j3) -- (A\i\j0);
}}
\foreach \i in {0,1,2,3}{
\pgfmathtruncatemacro{\ii}{\i+1};
\foreach \j in {0,1,2,3,4}{
\draw[thick] (A\i\j3) -- (A\ii\j1);
}}
\foreach \i in {0,1,2,3,4}{
\foreach \j in {0,1,2,3}{
\pgfmathtruncatemacro{\jj}{\j+1};
\draw[thick] (A\i\j0) -- (A\i\jj2);
}}
\foreach \i in {0,1,2,3,4}{
\draw[thick] (A\i02) to[bend left = 21] (A\i40);
}
\foreach \j in {0,1,2,3,4}{
\draw[thick] (A0\j1) to[bend left = 20] (A4\j3);
}
\end{scope}
\end{tikzpicture}
}
\caption{A different presentation of the graph $\cXb(5,20,1,4,18)$.}
\label{fig:example_squares}
\end{center}
\end{figure}
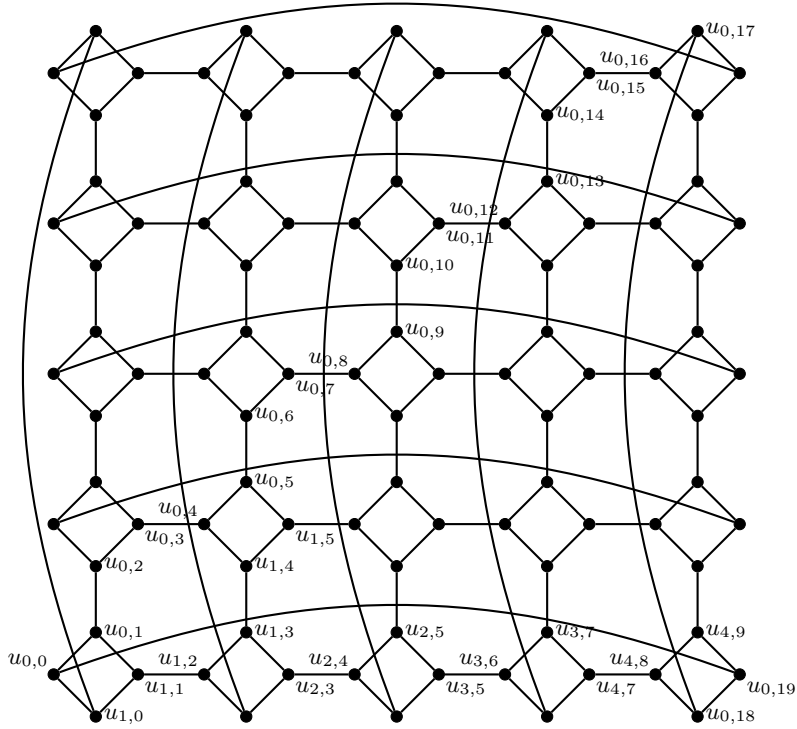

It can be verified that for each odd $m \geq 3$ the graph $\cXb(m,4m,1,4,4m-2)$ does indeed admit an automorphism $\eta$ with the above properties. Instead of exhibiting the explicit action of $\eta$ and proving that it is indeed an automorphism of the graph, which is a rather tedious task, we present an alternative presentation of the graph $\cXb(m,4m,1,4,4m-2)$ for $m = 5$ (see Figure~\ref{fig:example_squares}), from which the existence of $\eta$ is clear. It is also clear that a presentation of $\cXb(m,4m,1,4,4m-2)$ of this kind can be obtained for each odd $m \geq 3$. This finally gives the following result.

\begin{proposition}
\label{pro:m_odd_l_even}
Suppose $m$, $n$, $a$, $b$ and $\ell$ satisfy the assumptions and conditions of item (iii) of Theorem~\ref{the:main}, and suppose in addition that $m$ is odd and $\ell$ is even. Then the automorphism group of the graph $\G = \cXb(m,n,a,b,\ell)$ is regular and the $2$-factor $\cC$ is $\Aut(\G)$-invariant, unless $\G = \cXb(m,4m,1,4,4m-2)$, in which case $\G$ is vertex-transitive with vertex-stabilizers of order $2$ and $\cC$ is not $\Aut(\G)$-invariant.
\end{proposition}

We conclude the paper with two suggestions for possible future research. The first task is of course to determine the automorphism groups of the remaining two subfamilies of the graphs $\cXb(m,n,a,b,\ell)$ from item (iii) of Theorem~\ref{the:main}.

\begin{problem}
Determine which of the graphs from item (iii) of Theorem~\ref{the:main} have regular automorphism groups and determine the full automorphism groups of those that do not.
\end{problem}

Finally, we suggest a natural possibility of how to extend the work of~\cite{AlsKreKho19, AlsDobKhoSpa22, AlsSpa24} and the present paper in the efforts towards new results on~\cite[Problem~1.1]{AlsSpa24}. Namely, in~\cite{AlsSpa24} and the present paper we were concerned with the possibility that the quotient graph $\G_\cC$ of the given cubic graph $\G$ with respect to the given $2$-factor $\cC$ is a cycle. Another natural possibility to consider is when $\G_\cC$ is a complete graph. Of course, the case that $\G_\cC$ is $K_2$ is precisely the situation considered in~\cite{AlsDobKhoSpa22}, while the case that it is $K_3$ is covered by~\cite{AlsSpa24} and the present paper. One can thus restrict to the possibility that the quotient graph is of order at least $4$.

\begin{problem}
Classify the cubic vertex-transitive graphs $\G$, admitting a $2$-factor $\cC$ with at least $4$ components and a vertex-transitive group $G$ preserving $\cC$, such that the corresponding quotient graph $\G_\cC$ is a complete graph. For all such graphs determine whether the corresponding $2$-factor is $\Aut(\G)$-invariant and determine their full automorphism groups.
\end{problem}

\section*{Declarations and Acknowledgments}

\noindent
{\bf Competing interests:} The author has no relevant competing interests to declare.
\bigskip

\noindent
{\bf Data availability:} There is no associated data.
\bigskip

\noindent
{\bf Funding:} The author acknowledges the financial support by the Slovenian Research and Innovation Agency (research program P1-0285 and research project J1-50000).


\begin{thebibliography}{}
\bibitem{Als21} B.~Alspach, 
	Honeycomb toroidal graphs, 
	{\em Bull. Inst. Comb. Appl.} {\bf 91} (2021), 94--114.
\bibitem{AlsDea09} B.~Alspach, M.~Dean,
	Honeycomb toroidal graphs are Cayley graphs,
	{\em Inform. Process. Lett.} {\bf 109} (2009), 705--708.
\bibitem{AlsKreKho19} B.~Alspach, A.~Khodadadpour, D.L.~Kreher,
	On factor-invariant graphs,
	{\em Discrete Math.} {\bf 342} (2019), 2173--2178.
\bibitem{AlsDobKhoSpa22} B.~Alspach, T.~Dobson, A.~Khodadadpour, P.~\v Sparl,
	On factor-invariant graphs with two cycles,
	{\em Discrete Math.} {\bf 345} (2022), Paper No. 112937, 13 pp.
\bibitem{AlsSpa24} B.~Alspach, P.~\v Sparl,
	Cubic factor-invariant graphs of cycle quotient type -- The alternating case,
	{\em European J. Combin.} {\bf 120} (2024), Paper No. 103964, 22 pp.
\bibitem{BarGraSpi25} M.~Barbieri, V.~Grazian, P.~Spiga,	
	On the order of semiregular automorphisms of cubic vertex-transitive graphs,
	{\em European J. Combin.} {\bf 124} (2025), Paper No. 104091, 20 pp.
\bibitem{DobHujImrOrt25} T.~Dobson, A.~Hujdurovi\'c, W.~Imrich, R.~Ortner,
	On cubic vertex-transitive graphs of given girth,
	{\em Art Discrete Appl. Math.} {\bf 8} (2025), no. 1, Paper No. 1.12, 16 pp.
\bibitem{EibJajSpa19} E.~Eiben, R.~Jajcay, P.~\v Sparl,
	Symmetry properties of generalized graph truncations,
	{\em J. Combin. Theory, Ser. B} {\bf 137} (2019), 291--315.
\bibitem{GuoMoh24} K.~Guo, B.~Mohar,
	Simple eigenvalues of cubic vertex-transitive graphs,
	{\em Canad. J. Math.} {\bf 76} (2024), no. 5, 1496--1519.
\bibitem{KutMarMikSpa24} K.~Kutnar, D.~Maru\v si\v c, \v S.~Miklavi\v c, P. \v Sparl,
	On the structure of consistent cycles in cubic symmetric graphs,
	{\em J. Graph Theory} {\bf 105} (2024), no. 3, 337--356.
\bibitem{LiKwoZho23} N.~Li, Y.~S.~Kwon, J.-X.~Zhou,
	On cubic vertex-transitive non-Cayley graphs of $2$-power orders,
	{\em Discrete Math.} {\bf 346} (2023), no. 12, Paper No. 113625, 12 pp.
\bibitem{MikPotWil08} \v S.~Miklavi\v c, P.~Poto\v cnik, S.~Wilson,
	Arc-transitive cycle decompositions of tetravalent graphs,
	{\em J. Combin. Theory Ser. B} {\bf 98} (2008), 1181--1192
\bibitem{PotSpiVer13} P.~Poto\v cnik, P.~Spiga, G.~Verret, 
	Cubic vertex-transitive graphs on up to 1280 vertices, 
	{\em J. Symbolic Comput.} {\bf 50} (2013), 465--477.
\bibitem{PotTol23} P.~Poto\v cnik, M.~Toledo,
	Cubic vertex-transitive graphs admitting automorphisms of large order,
	{\em Bull. Malays. Math. Sci. Soc.} {\bf 46} (2023), no. 4, Paper No. 133, 33 pp.
\bibitem{PotWil07} P.~Poto\v cnik, S.~Wilson,
	Tetravalent edge-transitive graphs of girth at most 4,
	{\em J. Combin. Theory, Ser. B} {\bf 97} (2007), 217--236.
\bibitem{SirTucWat01} J.~\v Sir\'an, T.~W.~Tucker, M.~E.~Watkins, 
	Realizing finite edge-transitive orientable maps, 
	{\em J. Graph Theory} {\bf 37} (2001), 1--34.
\bibitem{Spa22} P.~\v Sparl,
	Symmetries of the honeycomb toroidal graphs,
	{\em J. Graph Theory} {\bf 99} (2022), 414--424.
\end{thebibliography}
\end{document}